\documentclass[12pt]{amsart}

\usepackage{subfiles}
\usepackage{mathrsfs}

\usepackage{hyperref}

\usepackage{times}
\usepackage{amssymb}

\usepackage{enumitem}

\usepackage{amsthm}
\usepackage{thmtools,thm-restate}
\usepackage[capitalise,noabbrev]{cleveref}

\declaretheorem[numberwithin=section]{lemma}
\declaretheorem[sibling=lemma]{theorem}

\declaretheorem[sibling=lemma]{corollary}
\declaretheorem[sibling=lemma]{proposition}

\declaretheorem[sibling=lemma,style=remark]{example}

\declaretheorem[sibling=lemma,style=remark]{definition}

\Crefname{inconvenientdefinition}{inconvenientdefinition}{{inconvenientdefinition}}

\usepackage{color}
\usepackage{gensymb}
\usepackage[usenames,dvipsnames,svgnames,table]{xcolor}

\usepackage{tikz}
\usepackage[backgroundcolor=yellow]{todonotes}
\usetikzlibrary{matrix,arrows,calc}
\usetikzlibrary{decorations.markings}
\usetikzlibrary{patterns}
\usetikzlibrary{snakes}

\usepackage{subfig}
\usepackage{xspace}
\usepackage[T1]{fontenc}
\usepackage[utf8]{inputenc}
\usepackage{bm}

\usepackage{subfiles}

\newcommand{\bigexists}{%
\mathop{\lower0.75ex\hbox{%
   \scalebox{1.7}{\ensuremath{\exists}}}}\limits}

\DeclareMathOperator{\St}{St}

\newcommand{\Sym}[1]{\mathfrak{S}(#1)}

\newcommand{\Embed}{\mathfrak{N}}

\DeclareMathOperator{\genus}{genus}

\newcommand{\V}{V}
\newcommand{\WC}{\mathcal{W}}
\newcommand{\BC}{\mathcal{B}}
\newcommand{\loops}{\mathcal{L}}
\newcommand{\E}{\mathcal{E}}
\newcommand{\F}{\mathcal{F}}

\newcommand{\Q}{\mathbb{Q}}

\newcommand{\AllYoung}{\mathbb{Y}}

\DeclareMathOperator{\Ch}{Ch}
\newcommand{\Chtt}{\Ch^{\ttwisted}}
\DeclareMathOperator{\ttwisted}{top}
\DeclareMathOperator{\ttop}{top}

\DeclareMathOperator{\twist}{twist}

\DeclareMathOperator{\weight}{mon}

\newcommand{\Laurent}{\Q\left[A,A^{-1}\right]}

\numberwithin{equation}{section}

\newcommand{\faceA}{red!50}
\newcommand{\faceB}{blue}
\newcommand{\faceAfill}{red!10}
\newcommand{\faceBfill}{blue!20}

\thanks{Version identifier: \texttt{e9e541d
}}

\author{Agnieszka Czyżewska-Jankowska}
\address{Piławska 20/16, 50-538 Wrocław, Poland}
\email{agnieszka.czyzewska@gmail.com}

\author {Piotr \'Sniady}
\address{
Institute of Mathematics, Polish Academy of Sciences,
\mbox{ul.~\'Sniadec\-kich 8,} 00-956 Warszawa, Poland} 
\email{psniady@impan.pl}

\title[Bijection between oriented and non-oriented maps]%
{
Bijection between oriented maps \\ and weighted non-oriented maps}

\begin{document}

\begin{abstract}
We consider bicolored maps, i.e.~graphs which are drawn on surfaces,
and construct a bijection between 
(i) \emph{oriented} maps with \emph{arbitary} face structure,
and (ii) 
(weighted) \emph{non-oriented} maps with \emph{exactly one face}.
Above, each non-oriented map is counted with a multiplicity which is
based on the concept of \emph{the orientability generating series} and 
\emph{the measure of orientability of a map}.
This bijection has the remarkable property of preserving the underlying bicolored graph.
Our bijection shows equivalence between two explicit formulas for the top-degree
of Jack characters, i.e.~(suitably normalized) coefficients in the expansion of Jack symmetric functions 
in the basis of power-sum symmetric functions. 
\end{abstract}



\maketitle

\newcommand{\G}{G}

\setcounter{section}{-1}

\newcommand{\kolorSigma}{blue}
\newcommand{\kolorPi}{red}
\newcommand{\kolorV}{DarkGreen}
\newcommand{\kolorW}{RawSienna}

\section{Prologue}
\label{sec:prologue}

In order to motivate the Reader and to give her some flavor of the results to expect,
we shall present now some selected highlights
before getting involved in somewhat lengthy definitions.
We also deliberately postpone the bibliographic details.

\subsection{Maps and orientability}
\label{sec:intro-maps}

Roughly speaking, a \emph{map} $M=(\G,S)$ is a bicolored graph $\G$ which is drawn on a surface $S$.
We require that the connected components of $S\setminus \G$, called \emph{faces}, 
are all homeomorphic to open discs.
The set of vertices of $\G$ is decomposed into two disjoint sets: the set of white vertices
and the set of black vertices. Each edge connects two vertices of opposite colors;
multiple edges are allowed. We do not allow isolated vertices.
We allow the surface $S$ to be disconnected.

\medskip

\begin{figure}
\centerline{
\begin{tikzpicture}[scale=0.6,
white/.style={circle,thick,draw=black,fill=white,inner sep=4pt},
black/.style={circle,draw=black,fill=black,inner sep=4pt},
connection/.style={draw=black,thick,black,auto}
]
\scriptsize
\begin{scope}
\clip (0,0) rectangle (10,10);
\fill[pattern color=\faceAfill,pattern=north west lines] (0,0) rectangle (10,10);
\fill[\faceBfill] (5,5) rectangle (8,8);
\draw (2,7)  node (b1)     [black] {};
\draw (3,2)  node (w1)     [white] {};
\coordinate (w1prim)      at (13,2);
\coordinate (w1bis)       at (3,12);
\coordinate (b1bis)       at (2,-3);
\draw (7,3)  node (bb1)    [black] {};
\coordinate (bb1prim) at (-3,3);
\draw (5,5)  node (AA)     [black] {};
\draw (8,5)  node (BA)     [white] {};
\draw (5,8)  node (AB)     [white] {};
\draw (8,8)  node (BB)     [black] {};
\draw[connection]         (w1)      to node {$5$}  (AA);
\draw[connection]         (AA)      to node {$3$}  (AB);
\draw[connection]         (AB)      to node {$8$}  (BB);
\draw[connection]         (BB)      to node {$6$}  (BA);
\draw[connection]         (BA)      to node [swap] {$2$} (AA);
\draw[connection]         (b1)       to  node [swap,pos=0.3] {$7$} (w1);
\draw[connection]         (bb1)     to node [swap,pos=0.3] {$1$} (w1prim);
\draw[connection]         (w1)      to node [pos=0.2]      {$1$} (bb1prim);
\draw[connection]         (w1)      to node[pos=0.15] {$4$}  (b1bis);
\draw[connection]         (b1)      to node [swap]    {$4$} (w1bis);
\draw[connection]         (w1)      to node [swap]    {$9$} (bb1);
\end{scope}
\draw[very thick,decoration={
    markings,
    mark=at position 0.666  with {\arrow{>}}},
    postaction={decorate}]  
(0,0) -- (10,0);
\draw[very thick,decoration={
    markings,
    mark=at position 0.666  with {\arrow{>}}},
    postaction={decorate}]  
(0,10) -- (10,10);
\draw[very thick,decoration={
    markings,
    mark=at position 0.666  with {\arrow{>>}}},
    postaction={decorate}]  
(10,0) -- (10,10);
\draw[very thick,decoration={
    markings,
    mark=at position 0.666  with {\arrow{>>}}},
    postaction={decorate}]  
(0,0) -- (0,10);
\end{tikzpicture}
}
\caption{Example of an \emph{oriented} map drawn on the torus.
The left side of the square should be glued  to the right side, 
as well as bottom to top, as indicated by the arrows.
The statement that it is an \emph{oriented} map means that
there is a prescribed consistent choice around each vertex
of what \emph{clockwise} means.}
\label{fig:torus}
\end{figure}
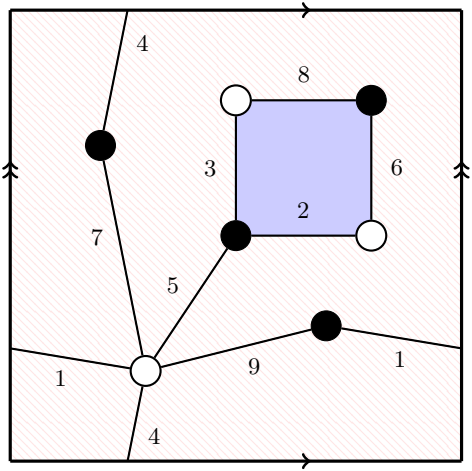

\begin{figure}
\centerline{
\begin{tikzpicture}[scale=0.6,
black/.style={circle,thick,draw=black,fill=white,inner sep=4pt},
white/.style={circle,draw=black,fill=black,inner sep=4pt},
connection/.style={draw=black,black,auto}
]
\scriptsize
\begin{scope}
\clip (0,0) rectangle (10,10);
\fill[pattern color=\faceAfill,pattern=north west lines] (0,0) rectangle (10,10);
\fill[\faceBfill] (5,5) rectangle (8,8);
\draw (2,7)  node (b1)     [black] {};
\draw (12,3) node (b1prim) [black] {};
\draw (3,2)  node (w1) [white] {};
\draw (-7,8) node (w1prim) [white] {};
\draw (5,5)  node (AA)     [black] {};
\draw (8,5)  node (BA)     [white] {};
\draw (5,8)  node (AB)     [white] {};
\draw (8,8)  node (BB)     [black] {};
\draw[connection]         (w1)      to node [pos=0.6] {$4$} node [swap,pos=0.6] {$9$} (AA);
\draw[connection]         (AA)      to node {$5$} node [swap] {$D$} (AB);
\draw[connection]         (AB)      to node {$6$} node [swap] {$C$} (BB);
\draw[connection]         (BB)      to node {$7$} node [swap] {$B$} (BA);
\draw[connection]         (BA)      to node {$8$} node [swap] {$A$} (AA);
\draw[connection]         (w1)      to node [pos=0.5] {$10$} node [swap,pos=0.425] {$2$} (b1prim);
\draw[connection]         (w1prim)  to node [pos=0.87] {$2$} node [swap,pos=0.95] {$10$} (b1);
\draw[connection]         (w1)      to node [pos=0.625] {$1$}  node [swap,pos=0.5] {$3$} (b1);
\end{scope}
\draw[very thick,decoration={
    markings,
    mark=at position 0.666  with {\arrow{>}}},
    postaction={decorate}]  
(0,0) -- (10,0);
\draw[very thick,decoration={
    markings,
    mark=at position 0.666  with {\arrow{>}}},
    postaction={decorate}]  
(10,10) -- (0,10);
\draw[very thick,decoration={
    markings,
    mark=at position 0.666  with {\arrow{>>}}},
    postaction={decorate}]  
(10,0) -- (10,10);
\draw[very thick,decoration={
    markings,
    mark=at position 0.666  with {\arrow{>>}}},
    postaction={decorate}]  
(0,10) -- (0,0);
\end{tikzpicture}}
\caption{Example of a \emph{non-oriented map} drawn on the projective plane.
The left side of the square should be glued with a twist to the right side, as well as bottom to top, 
as indicated by the arrows. The two faces of the map are indicated by the colors.}
\label{fig:map-projective-plane}

\centering
\begin{tikzpicture}[scale=1,
black/.style={circle,thick,draw=black,fill=white,inner sep=7pt},
white/.style={circle,draw=black,fill=black,inner sep=7pt},
connection/.style={auto,double distance=3mm}
]
\scriptsize

\begin{scope}

\draw[double distance=3mm,black!80]
 (5,5) -- 
node[right=3mm] {\textcolor{black}{$D$}} node [swap,left=3mm] {\textcolor{black}{$5$}} (5,8) --
node[below=3mm] {\textcolor{black}{$C$}} node [swap,above=3mm] {\textcolor{black}{$6$}} (8,8) --
node[right=3mm] {\textcolor{black}{$7$}} node [swap,left=3mm] {\textcolor{black}{$B$}} (8,5) --
node[below=3mm] {\textcolor{black}{$8$}} node [swap,above=3mm] {\textcolor{black}{$A$}} (5,5) --
node[below right=3mm] {\textcolor{black}{$9$}} node [swap,above left=3mm] {\textcolor{black}{$4$}} (3,3); 

\draw (0.5,2.5) node {$3$};
\draw (1.3,1.7) node {$1$};

\draw (2.5,0.5) node {$10$};
\draw (1.7,1.3) node {$2$};

\draw[black!80] plot[smooth] file {ribbonB.txt};
\draw[draw=white,line width=2mm] plot[smooth] file {ribbonA.txt};
\draw[black!80] plot[smooth] file {ribbonA.txt};

\draw[black!80] plot[smooth] file {ribbonC.txt};
\draw[black!80] plot[smooth] file {ribbonD.txt};

\draw (3,3)  node (w1) [white] {};
\draw (0,0)  node (w1) [black] {};

\draw (5,5)  node (AA)     [black] {};
\draw (8,5)  node (BA)     [white] {};
\draw (5,8)  node (AB)     [white] {};
\draw (8,8)  node (BB)     [black] {};

\end{scope}

\end{tikzpicture}
\caption{Alternative graphical representation of the map from \protect\cref{fig:map-projective-plane} as a ribbon graph. 
The edges of the map are represented as thin ribbons attached to the vertices.}
\label{fig:ribbon-graph}
\end{figure}

The maps which we consider in the current paper come in the following two flavors:
\emph{oriented maps} (which are drawn on an orientable surface $S$ which comes with
some prescribed choice of the orientation, see \cref{fig:torus}) 
and \emph{non-oriented maps} (which are drawn on an arbitrary surface $S$ without any 
additional structures, see \cref{fig:map-projective-plane}). 
These two flavors are quite distinct, nevertheless
some algebraic-combinatorial conjectures which concern \emph{Jack polynomials}
suggest that there is some hypothetical
natural one-parameter interpolation between them. 
To be more specific: it has been conjectured
that there exists some \emph{measure of non-orientability} of a given non-oriented map $M$
which is a hypothetical polynomial $\operatorname{weight}_M(\gamma)$ in the deformation parameter $\gamma$
with the property that
if a summation over \emph{non-oriented} maps is performed,
and each map is counted with appropriate multiplicity $\operatorname{weight}_M(\gamma)$,
the resulting sum becomes a natural interpolation between 
some generating series of oriented and non-oriented maps.

In the current paper we consider a concrete formula for such a candidate weight,
a candidate which will be denoted by $\weight_M(\gamma)$, which is an acronym that 
stands for \emph{\textbf{m}easure \textbf{o}f \textbf{n}on-orientability}.

\subsection{Top-degree of the measure of non-orientability}
\label{sec:mon-first}

In order to state the main results of the paper we will not need the definition
of the full polynomial $\weight_M(\gamma)\in\Q[\gamma]$ 
and we shall restrict ourselves to \emph{the leading coefficient} of this polynomial
which will be denoted by $\weight^{\ttop}_{M}$. We shall present now its definition.

\medskip

Any non-oriented map $M$ can be equivalently viewed as a \emph{ribbon graph},
see \cref{fig:ribbon-graph}.
With this viewpoint each vertex becomes a disc, 
each edge becomes a thin ribbon connecting the discs,
and each face of the map becomes a connected component of the boundary of the union of the discs
and the ribbons.

From a ribbon graph corresponding to $M$ we shall
remove all of its edges, one after another, in a
uniformly random order. 
We use the convention that whenever after an edge removal some vertex 
becomes isolated, we remove this vertex as well.

The following can be viewed as a \emph{working definition} of the quantity ${\weight^{\ttop}_{M}}$.

\begin{definition}
\label{def:mon-top-degree}
Assume that $M$ is a non-oriented map.
The quantity $\weight^{\ttop}_{M}$
is defined to be the probability of 
\emph{the event that in the above process of uniformly random edge removal,
at each step the number of faces of the ribbon graph is equal to the number of its connected 
components}.
\end{definition}
The link between the above quantity $\weight^{\ttop}_M$
and the polynomial $\weight_M$ will be provided later in \cref{prop:why-montop}.

\begin{figure}[t]
\centering
\begin{tikzpicture}[scale=1,
white/.style={circle,thick,draw=black,fill=white,inner sep=4pt},
black/.style={circle,draw=black,fill=black,inner sep=4pt},
]
\begin{scope}
\clip (0,0) rectangle (5,5);
\fill[blue!3] (0,0) rectangle (5,5);
\node (v1) at (1.3,3.5) [black] {};
\node (v2) at (3.7,1.5) [white] {};
\draw (v2) +(-5,-0) to node [auto,swap,pos=0.7] {\footnotesize $2$} node [auto,pos=0.8] {\footnotesize $4$}
        (v1) to node [auto,swap] {\footnotesize $1$} node [auto] {\footnotesize $5$}
        (v2);
\draw (v1) +(5,0) to node [auto,swap,pos=0.7] {\footnotesize $4$} node [auto,pos=0.8] {\footnotesize $2$}
        (v2);
\draw (v1) to node [auto] {\footnotesize $3$} node [auto,swap] {\footnotesize $6$} (1.3,5);
\draw (v2) to node [auto] {\footnotesize $3$} node [auto,swap] {\footnotesize $6$} (3.7,0);
\end{scope}
\draw[very thick,decoration={
    markings,
    mark=at position 0.666  with {\arrow{>}}},
    postaction={decorate}]  
(0,0) -- (5,0);

\draw[very thick,decoration={
    markings,
    mark=at position 0.666  with {\arrow{>}}},
    postaction={decorate}]  
(5,5) -- (0,5);

\draw[very thick,decoration={
    markings,
    mark=at position 0.666  with {\arrow{>>}}},
    postaction={decorate}]  
(0,0) -- (0,5);

\draw[very thick,decoration={
    markings,
    mark=at position 0.666  with {\arrow{>>}}},
    postaction={decorate}]  
(5,0) -- (5,5);

\end{tikzpicture}
\caption{
Example of a non-oriented map drawn on the Klein bottle: 
the left-hand side of the square should be glued to the right-hand side (without a twist) 
and the top side should be glued to the bottom side (with a twist), as indicated by the arrows.
This map has one face.}
\label{fig:exampleA}

\centering
\begin{tikzpicture}[scale=0.5,
white/.style={circle,thick,draw=black,fill=white,inner sep=7pt},
black/.style={circle,draw=black,fill=black,inner sep=7pt},
connection/.style={auto,double distance=3mm}
]
\clip (-3,-8) rectangle (13,5);
\draw[black!80,connection] (0,0) .. controls (0,10) and (10,-20) .. (10,0); 
\draw[connection,draw=white,ultra thick] (0,0) .. controls (0,-20) and (10,10) .. (10,0); 
\draw[black!80,connection] (0,0) .. controls (0,-20) and (10,10) .. (10,0); 
\fill[fill=white,draw=white] (0,-3mm) rectangle (10,3mm);

\draw[draw=white,ultra thick] plot[smooth] file {ribbonStraightA.txt};
\draw[black!80] plot[smooth] file {ribbonStraightA.txt};
\draw[draw=white,ultra thick] plot[smooth] file {ribbonStraightB.txt};
\draw[black!80] plot[smooth] file {ribbonStraightB.txt};

\draw (0,0) node[black] {};
\draw (10,0) node[white] {};
\draw (4.5,0.3) node[above] {$3$};
\draw (4.5,-0.3) node[below] {$6$};

\draw (0,-4) node[anchor=east] {$1$};
\draw (0.5,-4) node[anchor=west] {$5$};

\draw (10,-4) node[anchor=west] {$4$};
\draw (9.5,-4) node[anchor=east] {$2$};
\end{tikzpicture}
\caption{The map from Figure~\ref{fig:exampleA} drawn as a \emph{ribbon graph}.}
\label{fig:exampleB}
\end{figure}

\begin{example}
\label{example:A}
We consider the map $M$ shown in Figure~\ref{fig:exampleA};
the corresponding ribbon graph is shown in Figure~\ref{fig:exampleB}.

Consider the case when the ribbon marked $\{3,6\}$ is removed first;
then the remaining two ribbons form an annulus which has
two faces and only one connected component, thus the event considered
in \cref{def:mon-top-degree} does not hold.  

Consider now the remaining two cases when either 
the ribbon $\{1,5\}$ or the ribbon $\{2,4\}$ is removed first;
then the remaining two ribbons form the Möbius strip which has one face and
one connected component. After yet another ribbon removal the unique remaining ribbon
forms a ribbon graph which again consists of one face and one connected component.
It follows that the event considered
in \cref{def:mon-top-degree} holds.

In this way we found that the probability considered in \cref{def:mon-top-degree}
is equal to $\frac{2}{3}$ and thus ${\weight^{\ttop}_{M}}=\frac{2}{3}$.
\end{example}

\subsection{Rooted maps}
The notion of a \emph{rooted map} takes a different form when we speak about
\emph{oriented maps} than in the case of \emph{non-oriented maps}.

More specifically, by a \emph{rooted oriented map} we mean an oriented map 
in which one \emph{edge} is decorated. For example, one can take the map from
\cref{fig:torus} and remove the labels of all edges, except for the edge
labeled by the symbol $1$, and declare that this \emph{edge} is decorated.

By a \emph{rooted non-oriented map} we mean a non-oriented map
in which one of the \emph{edge-sides} is decorated
(each edge consists of two edge-sides). For example, one can take the map from
\cref{fig:map-projective-plane} and remove the labels of all edge-sides, 
except for the edge-side labeled by $1$, and declare that this \emph{edge-side} is decorated.

\subsection{The first main result}
The following is one of the main results of the current paper
(for the other one see \cref{theo:top-degre-equal}).
It states that a summation over \emph{oriented maps} with \emph{arbitrary face structure}
is equivalent to a weighted summation over \emph{non-oriented maps}
with \emph{exactly one face}.

\begin{theorem}[The first main result]
\label{theo:main-bijection}
For all integers $n\geq 1$ the following formal linear combinations of 
\emph{bicolored graphs} are equal:
\begin{equation}
\label{eq:main-thm}
  \sum_{M_1=(\G_1,S_1)} \G_1  = 
\sum_{M_2=(\G_2,S_2)} {\weight^{\ttop}_{M_2}}\ \G_2,    
\end{equation}
where the sum on the left-hand side runs over 
\emph{oriented, unlabeled, rooted, connected maps $M_1$ with $n$ edges}
and the sum on the right-hand side runs over 
\emph{non-oriented, unlabeled, rooted maps $M_2$ with $n$ edges and one face}.
\end{theorem}

\textbf{It should be stressed that on the left-hand side we impose no restrictions 
on the face-type of the map $M_1$, in particular the number of the faces
is arbitrary.}

Equation \eqref{eq:main-thm} is an equality between formal linear combinations
of bicolored graphs. Equivalently, it can be viewed as the following statement:
for each bicolored graph $G$, 
the number of ways in which $G$ can be drawn on 
\begin{enumerate}[label=(A1)]
   \item 
\label{A1}
some oriented surface $S_1$ in such a way that $(G,S_1)$ becomes an \emph{oriented map}
\end{enumerate}
is equal to the (weighted) number of ways in which $G$ can be drawn on 
\begin{enumerate}[label=(A2)]
   \item
\label{A2}
some non-oriented surface $S_2$
in such a way that $(G,S_2)$ becomes a \emph{non-oriented map with exactly one face}.
\end{enumerate}

\medskip

Our proof of \cref{theo:main-bijection} will be bijective.
For a fixed bicolored graph~$G$ we will find a bijection between:
\begin{enumerate}[label=(B\arabic*)]
   \item 
\label{B1}
the set of \emph{non-oriented (but orientable) maps $(G,S_1)$}, together with a choice of an \emph{arbitrary linear order} on the set of edges of $G$, and 
   \item 
\label{B2}
the set of \emph{non-oriented maps $(G,S_2)$}, 
together with a choice of a \emph{linear order on the set of edges of $G$ with the property
that if the edges of a ribbon graph corresponding to the map $(G,S_2)$ are removed 
according to this linear order, then the condition from \cref{def:mon-top-degree}
holds true}, 
i.e.~at each step the number of faces of the ribbon graph is equal to the number of its connected 
components.
\end{enumerate}

\medskip

Roughly speaking, our bijection 
consists of a number of \emph{twists} of the ribbons, see \cref{fig:twist}.
The proof of \cref{theo:main-bijection} is postponed to \cref{sec:proof:thm:main}.

\begin{figure}[t]

\centering
\subfloat[][]{
\begin{tikzpicture}[scale=0.4,rotate=90]
\clip (-3,-4) rectangle (13,4);
\draw[draw=black,thick] (0,-0.5) rectangle (10,0.5);
\begin{scope}[rotate=120]
   \draw[thick] (0,-0.5cm) rectangle (10,0.5);
   \draw (2.5, 1.2) node {$m$};
   \draw (2.5,-1.2) node {$n$};
\end{scope}
\begin{scope}[rotate=-120]
   \draw[thick] (0,-0.5cm) rectangle (10,0.5);
   \draw (2.5, 1.2) node {$k$};
   \draw (2.5,-1.2) node {$l$};
\end{scope}

\begin{scope}[xshift=10cm]
    \begin{scope}[rotate=60]
      \draw[thick] (0,-0.5) rectangle (10,0.5);
      \draw (2.5, 1.2) node {$z$};
      \draw (2.5,-1.2) node {$y$};
    \end{scope}
    \begin{scope}[rotate=-60]
      \draw[thick] (0,-0.5) rectangle (10,0.5);
      \draw (2.5, 1.2) node {$x$};
      \draw (2.5,-1.2) node {$w$};

    \end{scope}
\end{scope}

\fill (0,0) circle (1); 
\draw[fill=white,ultra thick] (10,0) circle (1); 
\draw (3,1.2)  node{$a$};
\draw (3,-1.2) node {$b$};

\end{tikzpicture}
\label{subfig:before-twist}
}
\hspace{2cm}
\subfloat[][]{
\begin{tikzpicture}[scale=0.4,rotate=90]
\clip (-3,-4) rectangle (13,4);
\draw[draw=black,thick] plot[smooth] file {ribbonStraightAA.txt};
\draw[draw=white,line width=2mm] plot[smooth] file {ribbonStraightBB.txt};
\draw[draw=black,thick] plot[smooth] file {ribbonStraightBB.txt};
\begin{scope}[rotate=120]
   \draw[thick] (0,-0.5cm) rectangle (10,0.5);
   \draw (2.5, 1.2) node {$m$};
   \draw (2.5,-1.2) node {$n$};
\end{scope}
\begin{scope}[rotate=-120]
   \draw[thick] (0,-0.5cm) rectangle (10,0.5);
   \draw (2.5, 1.2) node {$k$};
   \draw (2.5,-1.2) node {$l$};
\end{scope}

\begin{scope}[xshift=10cm]
    \begin{scope}[rotate=60]
      \draw[thick] (0,-0.5) rectangle (10,0.5);
      \draw (2.5, 1.2) node {$z$};
      \draw (2.5,-1.2) node {$y$};
    \end{scope}
    \begin{scope}[rotate=-60]
      \draw[thick] (0,-0.5) rectangle (10,0.5);
      \draw (2.5, 1.2) node {$x$};
      \draw (2.5,-1.2) node {$w$};

    \end{scope}
\end{scope}

\fill (0,0) circle (1); 
\draw[fill=white,ultra thick] (10,0) circle (1); 
\draw (3,1.2)  node {$a$};
\draw (3,-1.2) node {$b$};

\end{tikzpicture}
\label{subfig:after-twist}
}

\caption{
Twist of a ribbon in a ribbon graph. 
The details of the notation will be explained in \cref{{sec:figure6isused}}.
\protect\subref{subfig:before-twist} A part of a map. 
The pair partition describing the structure of black vertices is $\BC=\big\{ \{b,k\}, \{l,m\}, \{n,a\}, \dots \big\}$,
the pair partition describing the white vertices is $\WC=\big\{ \{b,w\}, \{x,y\}, \{z,a\}, \dots \big\}$ 
and the pair-partition describing the edges is 
$\E=\big\{ \{a,b\}, \{k,l\}, \{m,n\},  \{w,x\}, \{y,z\}, \dots\big\}$.
\protect\subref{subfig:after-twist}
The outcome of a twist of the edge $\{a,b\}$. 
Only the partition describing the structure of white vertices has changed and is equal to 
$\WC'=\big\{ \{a,w\}, \{x,y\}, \{z,b\}, \dots \big\}$.}
\label{fig:twist}
\end{figure}

\medskip

One possible motivation for \cref{theo:main-bijection} is purely aesthetical.
However, there is also another important motivation related to the study of 
\emph{Jack polynomials} and \emph{Jack characters}
which will be explored in \cref{sec:introduction} and, in particular, in
\cref{sec:motivations1,sec:motivations2}.

\subsection{Overview}
The structure of this paper is twofold.
The first part (Sections~\ref{sec:nonoriented}--\ref{sec:proof:thm:main}) 
is self-contained and devoted to the proof of the first main result, \cref{theo:main-bijection}.
The second part (Sections \ref{sec:introduction}--\ref{sec:final-section})
presents the the motivations, the wider context, 
the bibliographic details and the applications to Jack characters,
in particular to the second main result, \cref{theo:top-degre-equal}.

\newpage

\section{Non-oriented maps}
\label{sec:nonoriented}
\label{sec:measure-of-non-orientability}

The notations presented in this chapter are based on the work of Dołęga, F\'eray and 
Śniady \cite{DolegaFeraySniady2013}.

\subsection{Non-oriented maps, informal viewpoint}
\label{subsec:non-oriented}

If we draw an edge of a non-oriented map $M=(G,S)$ with a fat pen
(or, alternatively, if we regard an edge of a corresponding ribbon graph)
its boundary consists of two \emph{edge-sides}. 
The maps which we consider in the current paper have labeled edge-sides, see \cref{fig:map-projective-plane}; 
in other words each edge carries two labels, one on each of its sides.
We also assume that each label is unique.

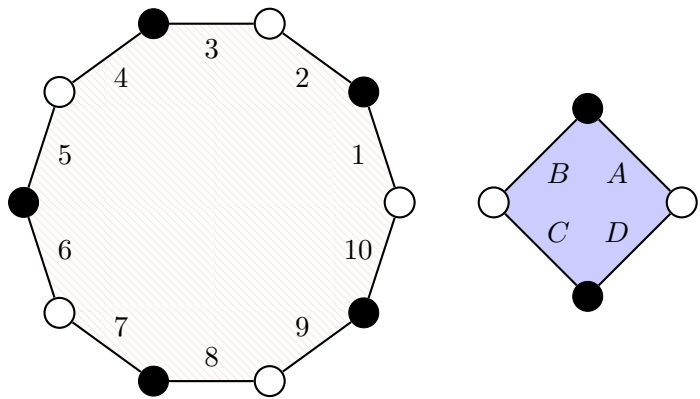
\begin{figure}[t]
\centering
\begin{tikzpicture}[scale=0.5,
black/.style={circle,thick,draw=black,fill=white,inner sep=4pt},
white/.style={circle,draw=black,fill=black,inner sep=4pt},
connection/.style={draw=black,black,auto}
]
\small

\fill[pattern color=\faceAfill,pattern=north west lines]  (0*36:5) -- (1*36:5) -- (2*36:5) -- (3*36:5) -- (4*36:5) -- (5*36:5) -- (6*36:5) -- (7*36:5) -- (8*36:5) -- (9*36:5);

\draw (0*36:5)  node (b1)     [black] {};
\draw (1*36:5)  node (b2)     [white] {};
\draw (2*36:5)  node (b3)     [black] {};
\draw (3*36:5)  node (b4)     [white] {};
\draw (4*36:5)  node (b5)     [black] {};
\draw (5*36:5)  node (b6)     [white] {};
\draw (6*36:5)  node (b7)     [black] {};
\draw (7*36:5)  node (b8)     [white] {};
\draw (8*36:5)  node (b9)     [black] {};
\draw (9*36:5)  node (b10)    [white] {};

\draw (0.5*36:4.1)  node {$1$};
\draw (1.5*36:4.1)  node {$2$};
\draw (2.5*36:4.1)  node {$3$};
\draw (3.5*36:4.1)  node {$4$};
\draw (4.5*36:4.1)  node {$5$};
\draw (5.5*36:4.1)  node {$6$};
\draw (6.5*36:4.1)  node {$7$};
\draw (7.5*36:4.1)  node {$8$};
\draw (8.5*36:4.1)  node {$9$};
\draw (9.5*36:4.1)  node {$10$};

\draw[connection] (b1) to  (b2);
\draw[connection] (b2) to  (b3);
\draw[connection] (b3) to  (b4);
\draw[connection] (b4) to  (b5);
\draw[connection] (b5) to  (b6);
\draw[connection] (b6) to  (b7);
\draw[connection] (b7) to  (b8);
\draw[connection] (b8) to  (b9);
\draw[connection] (b9) to  (b10);
\draw[connection] (b10)to (b1);

\begin{scope}[shift={(10,0)},scale=0.5]
\fill[\faceBfill]  (0:5) -- (90:5) -- (180:5) -- (270:5);

\draw (90*0:5)  node (b1)     [black] {};
\draw (90*1:5)  node (b2)     [white] {};
\draw (90*2:5)  node (b3)     [black] {};
\draw (90*3:5)  node (b4)     [white] {};

\draw (90*0.5:2.2) node {$A$};
\draw (90*1.5:2.2) node {$B$};
\draw (90*2.5:2.2) node {$C$};
\draw (90*3.5:2.2) node {$D$};

\draw[connection] (b1) to  (b2);
\draw[connection] (b2) to  (b3);
\draw[connection] (b3) to  (b4);
\draw[connection] (b4) to  (b1);
\end{scope}
\end{tikzpicture}

\caption{The map from \protect\cref{fig:map-projective-plane} is obtained from the above two polygons 
by gluing the following pairs of edges: 
$\{1,3\}$, $\{2,10\}$, $\{4,9\}$, $\{5,D\}$, $\{6,C\}$, $\{7,B\}$, $\{8,A\}$.
The colors of the polygons correspond to the colors of the faces of the map from 
\protect\cref{fig:map-projective-plane}.}
\label{fig:polygons}
\end{figure}

If we cut the surface $S$ along the edges of the graph $G$, the map becomes a collection of
bicolored polygons, each polygon corresponding to one face of the map, see \cref{fig:polygons}.
The original map $M$ can be recovered from this collection of polygons by gluing together
pairs of the edges of the polygons;
for each pair of glued edges  
the white (respectively, black) endpoints of the two edges should be glued together.

Thus a non-oriented map can be described alternatively as a collection of bicolored polygons with labeled edges 
(such as the ones from \cref{fig:polygons}) together with the pair-partition $\E$ telling which pairs of edges 
should be glued together.

\subsection{Non-oriented maps, more formal viewpoint}

\subsubsection{Pairings and polygons}
\label{SubsectPairingsPolygons}
A \emph{set-partition} of a set $X$ is a collection
$\{I_1,\dots,I_r\}$ of pairwise disjoint, non-empty subsets,
the union of which is equal to $X$.
A \emph{pairing} (or, alternatively, \emph{pair-partition}) of $X$ 
is a set-partition into pairs.

Let us consider now two pairings $\BC,\WC$ of the same set $X$ consisting of $2n$ elements.
We consider the following bicolored, edge-labeled graph $\loops(\BC,\WC)$:
\begin{itemize}
    \item it has $n$ black vertices indexed by the pairs of $\BC$
        and $n$ white vertices indexed by the pairs of $\WC$;
    \item its edges are labeled by the elements of $X$.
        The extremities of the edge labeled $i$ are the unique pair of $\BC$ containing $i$
        and the unique pair of $\WC$ containing $i$.
\end{itemize}

Note that each vertex has degree $2$ and each edge has one white and one black
extremity. Besides, if we erase the indices of the vertices, it is easy to
recover them from the labels of the edges (the index of a vertex is the set of
the two labels of the edges incident to this vertex). Thus, in the following we forget the indices
of the vertices and view $\loops(\BC,\WC)$ as an edge-labeled graph.

As every vertex has degree $2$, the graph $\loops(\BC,\WC)$ can be seen as a collection of
polygons.

\begin{example}
\label{example:1}
For partitions
\begin{align*}
    \BC &= \big\{ \{1,2\},\{3,4\},\{5,6\},\{7,8\},\{9,10\}, \{A,B\},\{C,D\} \big\}, \\
    \WC &= \big\{ \{2,3\},\{4,5\},\{6,7\},\{8,9\},\{10,1\}, \{B,C\},\{D,A\} \big\},
\end{align*}
the corresponding polygons $\loops(\BC,\WC)$ are shown in Figure~\ref{fig:polygons}.
\end{example}

\subsubsection{Non-oriented maps}

\label{subsec:non-oriented-appendix}

\begin{definition}
\label{def:set-theoretic-definitions}
    A \emph{non-oriented map} 
    is a triple $M=(\BC,\WC,\E)$ of pairings of the same base set $X$.
    \label{DefUmaps}
\end{definition}

The terminology comes from the fact that it is possible
to represent such a triple of pair-partitions as a bicolored graph embedded in a non-oriented 
(and possibly non-connected) surface.
Let us explain how this works.

\medskip

First, we consider the union of the polygons $\loops(\BC,\WC)$
defined above in \cref{SubsectPairingsPolygons}.
The edges of these polygons, that is the elements of the set~$X$,
are called \emph{edge-sides}.

We consider the union of the interiors of these polygons as a (possibly disconnected) surface
with a boundary.
If we consider two edge-sides, we can \emph{glue} them:
that means that we identify  with each other their white extremities, 
their black extremities, and the edge-sides themselves.

For any pair in the pairing $\E$, we glue the two corresponding edge-sides.
In this way we obtain a (possibly disconnected, possibly non-orientable)
surface $S$ without boundary.
After the gluing, the edges of the polygons form a bicolored graph $G$ embedded in the surface.
For instance, with the pairings $\BC$ and $\WC$ from \cref{example:1} and
\begin{equation}
\label{eq:example-pairing}
 \E = \big\{ \{1,3\},\{2,10\},\{4,9\},\{5,D\},\{6,C\}, \{7,B\},\{8,A\} \big\},  
\end{equation}
we get the graph from Figure~\ref{fig:map-projective-plane} embedded in the projective plane.

In general, the graph $G$ has as many connected components as the surface $S$.
Besides, the connected components of $S \setminus G$
correspond to the interiors of the collection of polygons we are starting from,
and, thus, they are homeomorphic to open discs.
These connected components are called \emph{faces}. 

This makes the link with the more common definition of maps:
usually, a (bicolored) map is defined as a (bicolored)
\emph{connected} graph $G$ embedded in a (non-oriented) surface $S$ in such a way that
each connected component of $S \setminus G$ is homeomorphic to an open disc.
It should be stressed that with our definition --- contrary to the traditional convention --- 
we do not require the map to be \emph{connected}.

Note that our maps have labeled edge-sides, and
each element of the label set $X$ is used exactly once as a label.

The pairing $\BC$ (respectively $\WC$) indicates which edge-sides share
the same corner around a black (respectively white) vertex.
This explains the names of these pairings.

\medskip

This encoding of non-oriented maps by triples of pairings
is of course not new.
It can for instance be found in \cite{Goulden1996a};
the presentation in that paper is nevertheless a bit different
as the authors consider there \emph{connected monochromatic} maps.

\subsection{Face-type. Summation over non-oriented maps}
\label{sec:summation-face-type}

Let $2\pi_1,2\pi_2,\dots,2\pi_\ell$ be the numbers of edges of the polygons, or --- equivalently --- 
the numbers of edges of the faces of the map. We say that $\pi=(\pi_1,\dots,\pi_\ell)$ 
is the \emph{face-type} of the collection of polygons or the \emph{face-type} of the map.

\medskip

The \emph{summation over non-oriented maps $M$ with face-type $(n)$} 
should be understood as follows: we fix a bicolored polygon $\loops$ with $2n$ labeled edges 
and consider all pair-partitions $\E$ of its edges; 
we sum over the resulting collection of maps $M=M(\E)$.
We will refer to this kind of summation as
\emph{conservative summation}.

\subsection{Edge liberation for non-oriented maps}
\label{sec:interesting-proof-begin}
Sometimes it will be convenient to consider a different way of summing over maps 
with the face-type $(n)$ in which the arrangement of the labels on the polygon $\loops$ is not fixed.
To be more specific: we consider all maps $M=(\BC',\WC',\E')$ where $\BC'$, $\WC'$ and $\E'$ are pair-partitions
of the same base set $X=[2n]$ with the property that $\loops(\BC',\WC')$ consists of a single polygon.
In order to differentiate this kind of summation, we will call it \emph{liberal summation} 
over non-oriented maps with the face-type~$(n)$.

\begin{proposition}
\label{prop:map-liberation-nonoriented}
For each integer $n\geq 1$ the following two formal linear combinations of
unlabeled maps are equal:
\begin{equation}
\label{eq:conservative-liberal}
\sum_{\substack{M \\ \text{liberal summation}}} \!\!\!\!\! M=
(2n-1)!  \!\!\!\!\!\!\!\!\!\!
\sum_{\substack{M \\ \text{conservative summation}}} \!\!\!\!\!\!\!\!\!\! M,   
\end{equation}
where on the left-hand side we consider a \emph{liberal summation} over non-oriented maps with the face-type $(n)$ 
and on the right-hand side we consider a \emph{conservative summation} over non-oriented maps with the face-type $(n)$.
On both sides of the equality we remove the labels of the edge-sides of the maps.
\end{proposition}
\begin{proof}
We fix the pair partitions $\BC$ and $\WC$ of the same set $X=[2n]$
in such a way that $\loops=\loops(\BC,\WC)$ is a single polygon. Let us consider the collection
\begin{equation}
\label{eq:collection}
\big\{ (\pi,\E) : \pi\in\Sym{2n} \text{ and $\E$ is a pair-partition of $X=[2n]$}\big\},
\end{equation}
where $\Sym{2n}$ denotes the symmetric group which we view as the set of permutations of $X$.
A formal sum of 
\begin{equation}
\label{eq:M}
M:=M(\E)   
\end{equation}
over this collection clearly corresponds to the conservative summation 
over non-oriented maps $M$ with face type $(n)$,
with each summand taken with the multiplicity $|\Sym{2n}|=(2n)!$;
thus this formal sum is equal to $2n$ times the right-hand side of \eqref{eq:conservative-liberal}.

\medskip

For a pair partition $I$ of the base set $X=[2n]$ and a permutation $\pi\in\Sym{2n}$
we denote by $\pi(I)$ the pair-partition of the same base set $X$ defined as follows:
for each $a,b\in S$ we have 
\[\{a,b\} \in I \iff \{\pi(a), \pi(b)\} \in \pi(I).\]

To each pair $(\pi,\E)$ we shall associate the map 
\begin{equation}
\label{eq:Mprim}
M':=(\BC',\WC',\E')
\end{equation}
with $\BC':=\pi(\BC)$, $\WC':=\pi(\WC)$, $\E':=\pi(\E)$.
The map $M'$ can be viewed as the map $M=(\BC,\WC,\E)$ with the permuted labels of the edge-sides;
in particular after removal of the labels of the edge-sides the maps $M$ and $M'$ are equal.
Thus the following two formal linear combinations of maps \emph{with removed labels} are equal:
\[  \sum M = \sum M',  \]
where both sums run over \eqref{eq:collection} while $M$, $M'$ should be understood as in
\eqref{eq:M}, \eqref{eq:Mprim}. 

\medskip

Clearly, the map $M'$ has one face and 
each non-oriented map $M'$ with one face on the base set $X$ can be obtained in this way.
Furthermore, two pairs
$(\pi,\E)$ and $(\sigma,\mathcal{F})$ have the same image $M'=(\BC',\WC',\E')$
if and only if $\mathcal{F}=\sigma^{-1} \pi(\E)$ and 
the permutation $\sigma^{-1}\pi$ leaves each of the partitions $\BC$ and $\WC$ invariant;
the set of such permutations forms the dihedral group $D_n$ of the isometries of a regular polygon with $n$ edges.
It follows that there is a bijective correspondence between the preimage of a given map $M'$
(we assume that $M'$ has a single face) and the dihedral group $D_n$.
Therefore a formal sum of $M':=(\BC',\WC',\E')$ over the collection \eqref{eq:collection}
corresponds to the liberal summation with each summand taken with the multiplicity $|D_n|=2n$;
in other words it is equal to the left-hand side of \eqref{eq:conservative-liberal}
multiplied by $2n$.

\smallskip

The comparison of the above two conclusions about the conservative and the liberal
ways of summing finishes the proof.
\end{proof}

\subsection{Removal of edges}
\label{sec:removal-of-edges}

If $E$ is an edge of a map $M$, we denote by $M\setminus E=M\setminus\{E\}$ 
\emph{the map $M$ with the edge $E$ removed}. 
This definition is a bit subtle; for example since we do not allow maps having isolated vertices, 
if some endpoint of $E$ is a leaf, we remove it as well. 
Removal of an edge might change the topology of the surface on which the map is drawn; 
for this reason instead of figures of the type presented in \cref{fig:map-projective-plane} 
it is more convenient to consider for this purpose \emph{ribbon graphs}, see \cref{fig:ribbon-graph}.

For a more rigorous treatment see \cite[Section 3.6]{DolegaFeraySniady2013}.

\subsection{Three kinds of edges}
\label{subsec:anatomy}

\begin{figure}[tbp]
\centering
\begin{tikzpicture}[scale=0.7,
black/.style={circle,thick,draw=black,fill=white,inner sep=4pt},
white/.style={circle,draw=black,fill=black,inner sep=4pt},
faceAs/.style={\faceA, dashed,  line width=6pt},
faceBs/.style ={\faceB, line width=6pt},
connection/.style={draw=black!80,black!80,auto}
]
\scriptsize

\begin{scope}
\clip (0,0) rectangle (10,10);

\fill[pattern color=\faceAfill,pattern=north west lines] (0,0) rectangle (10,10);

\draw (2,7)  node (b1)     [black] {};
\draw (12,3) node (b1prim) [black] {};
\draw (3,2)  node (w1) [white] {};
\draw (-7,8) node (w1prim) [white] {};

\draw (5,5)  node (AA)     [black] {};
\draw (8,5)  node (BA)     [white] {};
\draw (5,8)  node (AB)     [white] {};
\draw (8,8)  node (BB)     [black] {};

\draw[faceAs,left to-left to] (w1)      to (AA);
\draw[faceAs,left to-left to] (AA)      to (AB);
\draw[faceAs,left to-left to] (AB)      to (BB);
\draw[faceAs,left to-left to] (BB)      to (BA);
\draw[faceAs,left to-left to] (BA)      to (AA);
\draw[faceAs, ->] (w1)      to (b1prim);
\draw[faceAs, ->] (w1prim)  to (b1);
\draw[faceAs, <-] (w1)      to (b1);

\begin{scope}
\clip (5,5) -- (8,5) -- (8,8) -- (5,8);
\fill[\faceBfill] (5,5) rectangle (8,8);
\draw (5,5)  node (AA)     [black] {};
\draw (8,5)  node (BA)     [white] {};
\draw (5,8)  node (AB)     [white] {};
\draw (8,8)  node (BB)     [black] {};

\draw[faceBs, left to-left to] (AA) -- (AB);
\draw[faceBs, left to-left to] (AB) -- (BB);
\draw[faceBs, left to-left to] (BB) -- (BA);
\draw[faceBs, left to-left to] (BA) -- (AA);
\end{scope}

\draw[connection,draw=white,double=black,ultra thick]         (w1)      to node [pos=0.6] {\textcolor{black}{$4$}} node [swap,pos=0.6] {\textcolor{black}{$9$}} (AA);
\draw[connection,draw=white,double=black,ultra thick]         (AA)      to node {\textcolor{black}{$5$}} node [swap] {\textcolor{black}{$D$}} (AB);
\draw[connection,draw=white,double=black,ultra thick]         (AB)      to node {\textcolor{black}{$6$}} node [swap] {\textcolor{black}{$C$}} (BB);
\draw[connection,draw=white,double=black,ultra thick]         (BB)      to node {\textcolor{black}{$7$}} node [swap] {\textcolor{black}{$B$}} (BA);
\draw[connection,draw=white,double=black,ultra thick]         (BA)      to node {\textcolor{black}{$8$}} node [swap] {\textcolor{black}{$A$}} (AA);

\draw[connection,draw=white,double=black,ultra thick,pos=0.4] (w1)      to node [pos=0.5]  {\textcolor{black}{$10$}} node [swap,pos=0.425] {\textcolor{black}{$2$}} (b1prim);
\draw[connection,draw=white,double=black,ultra thick,pos=0.9] (w1prim)  to node [pos=0.82]  {\textcolor{black}{$2$}} node [swap,pos=0.9] {\textcolor{black}{$10$}} (b1);
\draw[connection,draw=white,double=black,ultra thick]         (w1)      to node [pos=0.625] {\textcolor{black}{$1$}} node [swap,pos=0.5] {\textcolor{black}{$3$}} (b1);

\end{scope}

\draw[very thick,decoration={
    markings,
    mark=at position 0.666  with {\arrow{>}}},
    postaction={decorate}]  
(0,0) -- (10,0);

\draw[very thick,decoration={
    markings,
    mark=at position 0.666  with {\arrow{>}}},
    postaction={decorate}]  
(10,10) -- (0,10);

\draw[very thick,decoration={
    markings,
    mark=at position 0.666  with {\arrow{>>}}},
    postaction={decorate}]  
(10,0) -- (10,10);

\draw[very thick,decoration={
    markings,
    mark=at position 0.666  with {\arrow{>>}}},
    postaction={decorate}]  
(0,10) -- (0,0);

\end{tikzpicture}
\caption{
The non-oriented map from \protect\cref{fig:map-projective-plane}.
On the boundary of each face some arbitrary orientation was chosen, as indicated by arrows.
The edge $\{4,9\}$ is an example of a \emph{straight edge}, 
the edge $\{1,3\}$ is an example of a \emph{twisted edge}, 
the edge $\{6,C\}$ is an example of an \emph{interface edge}.}
\label{fig:map-orientation}

\end{figure}
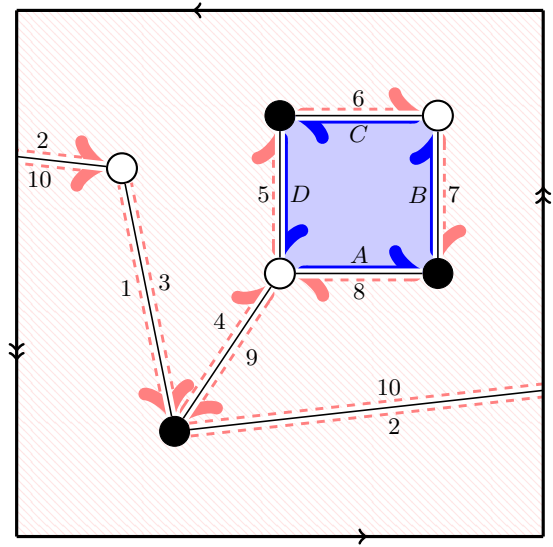

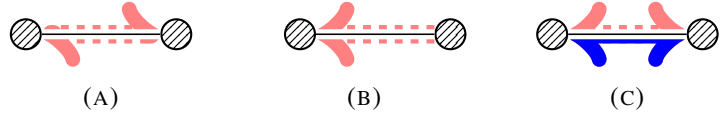
\begin{figure}[t]
\centering
\subfloat[]{
\begin{tikzpicture}[scale=0.5,
White/.style={circle,thick,draw=black,pattern=north east lines,inner sep=4pt},
Black/.style={circle,thick,draw=black,pattern=north east lines,inner sep=4pt},
connection/.style={draw=black!80,black!80,auto},
faceAs/.style={\faceA, dashed,  line width=6pt},
faceBs/.style ={\faceB, line width=6pt},
]

\draw[white] (-1,-1) rectangle (5,1);
\draw (0,0) node (b) [Black] {};
\draw (4,0) node (w) [White] {};

\draw[faceAs, line width=7pt, left to-left to] (w) to (b);
\draw[connection,draw=white,double=black,ultra thick] (w) to (b);
\end{tikzpicture}
\label{subfig:straight}
} 
\quad
\subfloat[]{
\begin{tikzpicture}[scale=0.5,
White/.style={circle,thick,draw=black,pattern=north east lines,inner sep=4pt},
Black/.style={circle,thick,draw=black,pattern=north east lines,inner sep=4pt},
connection/.style={draw=black!80,black!80,auto},
faceAs/.style={\faceA, dashed,  line width=6pt},
faceBs/.style ={\faceB, line width=6pt},
]

\draw[white] (-1,-1) rectangle (4,1);
\draw (0,0) node (b) [Black] {};
\draw (4,0) node (w) [White] {};

\draw[faceAs, line width=7pt, ->] (w) to (b);
\draw[connection,draw=white,double=black,ultra thick] (w) to (b);
\end{tikzpicture}
\label{subfig:mobius}
}
\quad
\subfloat[]{
\begin{tikzpicture}[scale=0.5,
White/.style={circle,thick,draw=black,pattern=north east lines,inner sep=4pt},
Black/.style={circle,thick,draw=black,pattern=north east lines,inner sep=4pt},
connection/.style={draw=black!80,black!80,auto},
faceAs/.style={\faceA, dashed,  line width=6pt},
faceBs/.style ={\faceB, line width=6pt},
]

\draw[white] (-1,-1) rectangle (4,1);
\draw (0,0) node (b) [Black] {};
\draw (4,0) node (w) [White] {};

\draw[faceAs, line width=7pt,<->] (w) to (b);
\begin{scope}
   \clip (-1,-1) rectangle (5,0);
   \draw[faceBs, line width=7pt, <->] (w) to (b);
\end{scope}

\draw[connection,draw=white,double=black,ultra thick] (w) to (b);
\end{tikzpicture}
\label{subfig:interface}
}
\caption{Three possible kinds of edges in a map (see \cref{fig:map-orientation}):
\protect\subref{subfig:straight} \emph{straight edge}: 
both edge-sides of the edge belong to the same face and have opposite orientations, 
\protect\subref{subfig:mobius} \emph{twisted edge}: 
both edge-sides of the edge belong to the same face and have the same orientation, 
\protect\subref{subfig:interface} \emph{interface edge}: 
the edge-sides of the edge belong to two different faces; their orientations are not important.
In all three cases the colors of the vertices are not important.}
\label{fig:3cases}
\end{figure}

Let a map $M$ with some selected edge $E$ be given (see example in \cref{fig:map-orientation}). 
There are three possibilities: 
\begin{itemize}
 \item \emph{Both sides of the edge $E$ are lying on the boundary of the same
face~$F$. This means that if we travel
along the boundary of the face $F$ then we visit the edge $E$ twice. Assume that
the directions in which we travel twice
along the edge $E$ are \emph{opposite}, see \cref{subfig:straight}. }

In this case the edge $E$ is called \emph{straight} and we associate to it the
weight $1$.

 \item \emph{Both sides of the edge $E$ are lying on the boundary of the same
face $F$. This means that if we travel
along the boundary of the face $F$ then we visit the edge $E$ twice. Assume that
the directions in which we travel twice
along the edge $E$ are \emph{the same}, see \cref{subfig:mobius}.}

In this case the edge $E$ is called \emph{twisted} and we associate to it the
weight $\gamma$.

 \item \emph{The edge $E$ is lying on the boundary of two different faces, see \cref{subfig:interface}.} 

In this case the edge $E$ is called \emph{interface} and we associate to it the
weight $\frac{1}{2}$. 
\end{itemize}

The weight given by the above convention will be denoted $\weight_{M,E}$;
we will need this notion much later in \cref{sec:weight-map-history}.
The classification of edges into three types (i.e.~\emph{straight} versus
\emph{twisted} versus \emph{interface}) will be necessary immediately.

\subsection{Top-degree histories}
\label{sec:histories}

The definitions presented in this section 
are directly related to the quantity $\weight^{\ttop}_M$ from \cref{def:mon-top-degree}.

\begin{definition}
\label{def:top-degree-map}
We say that a non-oriented map $M$ is a \emph{top-degree map} if each connected component of $M$ is one
of the faces of $M$.   
\end{definition}

For a given map $M$ we will say that a \emph{history} is an arbitrary 
linear order $\prec$ on the set of edges of $M$.
Let $E_1,\dots,E_n$ be the sequence of edges of $M$, listed according to the
linear order~$\prec$. We set 
$M_i = M \setminus \{ E_1,\dots, E_i \}$;
in other words $M_0,M_1,\dots,M_n$ is the sequence of maps obtained from $M$ by removing the edges,
one by one, in the order prescribed by the history $\prec$.

\begin{definition}
\label{def:top-degree}
We say that $(M,\prec)$ is a \emph{top-degree pair} if $M$ is a non-oriented map
and $\prec$ is a history with the property that
each of the maps $M_0,\dots,M_n$ defined above is a top-degree map.
\end{definition}

With the above definitions, for a given non-oriented map $M$ of face-type $(n)$,
the corresponding quantity $\weight^{\ttop}_M$ 
from \cref{def:mon-top-degree}
can be reformulated
as the probability that for a uniformly random choice of the history $\prec$,
the pair $(M,\prec)$ is a top-degree pair.

We will not make use of the following lemma
(but we will make use of its extension, \cref{lem:connected-component-one-face}).
Nevertheless, we decided to state it here because it provides a
natural and intuitive interpretation for the notion of \emph{top-degree pairs}.

\begin{lemma}
\label{lem:connected-component-one-face-abbridged}
Let $M$ be a map with $n$ edges and $\prec$ be a history.
We use the notations introduced above.
 
Then the following conditions are equivalent:
\begin{enumerate}[label={(\Alph*)}]

   \item \label{enum:maptopdegree-A}
   the pair $(M,\prec)$ is top-degree;

   \item \label{enum:types-of-edges-A}
      for each $0\leq i\leq n-1$ the edge $E_{i+1}$ of the map $M_i$ is either:
     \begin{itemize}
        \item a twisted edge, or,
        \item a bridge or a leaf.
     \end{itemize}

\end{enumerate}

\end{lemma}

Note that an alternative proof of this result (in a wider generality)
will be given in \cref{sec:top-degree-of-weight}.

\begin{proof}
Suppose that the condition \ref{enum:types-of-edges-A} holds true.
It is easy to prove by backward induction that
for each $0\leq i\leq n-1$ the map $M_i$ is a top-degree map
which implies the condition \ref{enum:maptopdegree-A}.

\medskip

Assume now that the condition \ref{enum:types-of-edges-A} does not hold true.
Let $i$ be the maximal number with the property that 
the edge $E_{i+1}$ of the map $M_i$ is neither
a twisted edge, nor a bridge, nor a leaf.
Clearly $i<n$.
The reasoning from the previous paragraph shows that the map $M_{i+1}$ 
is top-degree. The edge $E_{i+1}$ is not a bridge, so its endpoints must be located in 
the same connected component of $M_{i+1}$ therefore they belong to the same face.
Since the edge $E_{i+1}$ is not a twisted edge of $M_{i}$, it follows that 
$E_{i+1}$ is an interface edge of $M_{i}$. It follows that the connected component 
of the map $M_i$ which contains the edge $E_{i+1}$ consists of two faces.
In this way we proved that $M_i$ is not top-degree and the condition 
\ref{enum:maptopdegree-A} does not hold true.
\end{proof}

\subsection{Top degree histories and orientable maps}
\label{sec:figure6isused}

Let $M=(\BC,\WC,\E)$ be a non-oriented map and let $E$ be one of its edges.
We denote by $\twist_E(M)$ the non-oriented map which is obtained from $M$ by twisting the edge $E$,
see \cref{fig:twist}. Formally, $\twist_E(M)=(\BC,\WC',\E)$ is the map obtained by changing the structure
of the pairings of the white vertices, as explained in \cref{fig:twist}.

\begin{lemma}
\label{lem:dobra-przekrecalnosc}
Let $M$ be a map, let $E$ be one of its edges and let $M_1:=M\setminus E$.   
Assume that $M_1$ is a top-degree map (respectively, an orientable map).

If $E$ is a bridge or a leaf in the map $M$, then $M$
is a top-degree map (respectively, an orientable map).

If $E$ is neither a bridge nor a leaf in the map $M$, 
then exactly one of the following two maps: $M$ and $\twist_E(M)$ 
is a top-degree map (respectively, an orientable map). %Furthermore, the edge $E$ is a twisted edge or a straight edge
\end{lemma}

The proof is immediate.

\begin{theorem}
\label{theo:key-bijection}
Let an integer $n\geq 1$ be fixed and $X$ be an arbitrary label set with $|X|=2n$.
There exists a bijection $\Psi$ between:
\begin{enumerate}[label={(\alph*)}]
 \item \label{enum:twisted} the set of top-degree pairs $(M,\prec)$ where $M$ 
is a non-oriented map with $n$ edges and labels from $X$;

 \item \label{enum:orientable} the set of pairs $(M,\prec)$, 
where $M$ is an orientable, non-oriented map with $n$ edges and labels from $X$, 
and $\prec$ is a history on $M$.
\end{enumerate}
This bijection has the property that if $\Psi:(M,\prec)\mapsto (M',\prec')$ then the maps $M$ and $M'$ regarded as bicolored graphs are isomorphic.
\end{theorem}
\begin{proof}
The bijection which we will construct will have the following additional property:
for each $(M,\prec)$
there exists some subset $\{e_1,\dots,e_l\}$ of the set of the edges of $M$ 
such that
\[ \Phi(M,\prec)= \big( \twist_{e_1} \cdots \twist_{e_l} M, \prec \big).\]

\bigskip

We will use induction with respect to the variable $n$.
For $n=1$ both sets \ref{enum:twisted} and \ref{enum:orientable} consist of a single element
and there is an obvious bijection between them.

\bigskip

Consider the case $n\geq 2$. Let $(M,\prec)$ be a top-degree pair. Let $E$ be the first of the edges of $M$, 
according to the linear order $\prec$; let $M_1:=M\setminus E$ and let $\prec_1$ be the restriction of
the linear order $\prec$ to the edges of $M_1$; and let $X_1$ be the set $X$ with the labels of the edge-sides
of $E$ removed. 

The pair $(M_1,\prec_1)$ is also a top-degree pair and thus $\Phi(M_1,\prec_1)$ has been already 
constructed by the inductive assertion; 
there exists some set of the edges $\{e_1,\dots,e_l\}$ of $M_1$ 
with the property that
\[ \Phi(M_1,\prec_1)= \big( M_1', \prec_1 \big)\]
with
\[ M_1'= \twist_{e_1} \cdots \twist_{e_l} M_1.\]

By twisting the same set of edges in the bigger map $M$ we define
\[ \widetilde{M}:= \twist_{e_1} \cdots \twist_{e_l} M.\]

\medskip

There are the following two possibilities.
\begin{itemize}
\item \emph{The case when the edge $E$ is a bridge or a leaf in the map $M$.}

We define
\[ \Phi(M,\prec):= \big( \widetilde{M}, \prec \big).\]

Note that $M_1'=\widetilde{M}\setminus E$ is, by definition, an orientable map
and $E$ is a bridge or a leaf in the map $\widetilde{M}$; thus
\cref{lem:dobra-przekrecalnosc} implies that the map $\widetilde{M}$ is orientable, as required.

\item \emph{The case when the edge $E$ is neither a bridge nor a leaf in the map~$M$.}
      
Since $M_1'=\widetilde{M}\setminus E$ is an orientable map and $E$ 
is neither a bridge nor a leaf in the map $\widetilde{M}$,
\cref{lem:dobra-przekrecalnosc} implies that exactly one of the maps $\widetilde{M}$ and $\twist_E \widetilde{M}$
is orientable; we denote this orientable map by~$M'$.

Finally, we define
\[ \Phi(M,\prec):= \big( M', \prec \big).\]

   \end{itemize}

\medskip

This concludes the inductive construction.
It remains to show that the map $\Phi$ constructed above has an inverse.
The construction of the inverse $\Phi^{-1}$ is quite analogous to that
of $\Phi$ and we skip it.
\end{proof}

\section{Oriented maps}
\label{sec:oriented-maps-chapter}

\subsection{Oriented maps}
\label{sec:oriented-maps}
\label{sec:poczatek-zorientowanych-map}

We define an \emph{oriented} map as a bicolored
graph $G$ embedded in an \emph{oriented} surface $S$ in such a way that
each connected component of $S \setminus G$ is homeomorphic to an open disc.
If the number of edges is equal to $n$,
we shall assume that the edges are labeled by the elements of $[n]$ in such a way that
each label is used exactly once.

There is a bijective correspondence between such oriented maps with $n$ edges and
the set of pairs $(\sigma_1,\sigma_2)$, where $\sigma_1,\sigma_2\in\Sym{n}$ are permutations.
This correspondence follows from the observation that
the structure of such an oriented map is uniquely determined by the counterclockwise cyclic order of the edges 
around the white vertices 
(which we declare to be encoded by the disjoint cycle decomposition of the permutation~$\sigma_1$)
and by the counterclockwise cyclic order of the edges around the black vertices
(which we declare to be encoded by the disjoint cycle decomposition of the permutation~$\sigma_2$).
The corresponding oriented map will be denoted by $M(\sigma_1,\sigma_2)$.

\begin{example}
\label{example:map-on-torus}
The oriented map shown in \cref{fig:torus} corresponds to the pair 
\[ \sigma_1=(1,4,9,5,7)(2,6)(3,8), \qquad \sigma_2=(1,9)(2,3,5)(4,7)(6,8).\]
\end{example}

\subsection{Edge liberation for oriented maps}
\label{koniec-dowodu-o-skrecalnosci}

For $\sigma_1,\sigma_2\in\Sym{n}$ 
we say that \emph{``$\langle \sigma_1,\sigma_2 \rangle$ is transitive''} 
if the group generated by $\sigma_1$ and $\sigma_2$
acts transitively on the underlying set $[n]$.
It is easy to see that 
the oriented map $M(\sigma_1,\sigma_2)$ 
is \emph{connected} if and only if $\langle \sigma_1,\sigma_2 \rangle$ is transitive.

\medskip

The following result shows that a summation over
\emph{oriented, connected} maps 
can be alternatively viewed as a summation over \emph{orientable} (but not oriented), \emph{connected} maps.

\begin{proposition}
\label{prop:edge-liberation-oriented}
Let $n\geq 1$. 
Then the following two formal linear combinations of non-oriented unlabeled maps are equal:
\begin{equation}
\label{eq:liberal-conservative2}
(2n)!   \sum_{\substack{\sigma_1,\sigma_2\in \Sym{n} \\
\langle \sigma_1,\sigma_2 \rangle \text{ is transitive}}}  
M(\sigma_1,\sigma_2) = 
2 \cdot n!
\sum_M M,
\end{equation}
where the sum on the right-hand side runs over non-oriented maps
$M=(\BC,\WC,\E)$ where $\BC$, $\WC$, $\E$ are arbitrary pair-partitions 
of the set $X=[2n]$ such that $M$ is orientable and connected.
\end{proposition}
\begin{proof}
Consider the set    
\begin{multline}
\label{eq:transitive-pairs}
\big\{ (\sigma_1,\sigma_2, f) : \sigma_1,\sigma_2\in \Sym{n}\text{ and } 
\langle \sigma_1,\sigma_2 \rangle \text{ is transitive}, \\
f\colon[n]\times \{1,2\}\to [2n] \text{ is a bijection} \big\}.   
\end{multline}
A formal sum of $M(\sigma_1,\sigma_2)$ over this set 
is clearly equal to the left-hand side of \eqref{eq:liberal-conservative2}.

\medskip

Let us consider some triple which belongs to \eqref{eq:transitive-pairs}.
The permutations $\sigma_1$ and $\sigma_2$ define an \emph{oriented} map $M(\sigma_1,\sigma_2)$ 
with the edges labeled by the elements of $[n]$; 
in the following we will show how to view this map as a \emph{non-oriented} map.
This will be done by defining the labels associated to all edge-sides of the original map $M$;
these new labels belong to $[2n]$.
In other words, to each edge $e$ of $M$ (i.e., to each label $e\in[n]$)
we need to associate an ordered pair of two labels from $[2n]$; we declare this pair to be
$\big( f(e,1), f(e,2) \big)$. More specifically, if we go counterclockwise around the white endpoint of $e$
and read the labels of the edge-sides, then immediately after the label $f(e,1)$ we should read
the label $f(e,2)$, see \cref{fig:oriented-to-nonoriented}.

\newcommand{\BETA}{\sigma_1(\alpha)}
\newcommand{\GAMMA}{\sigma_1^{-1}(\alpha)}
\newcommand{\DELTA}{\sigma_2(\alpha)}
\newcommand{\EPSILON}{\sigma_2^{-1}(\alpha)}

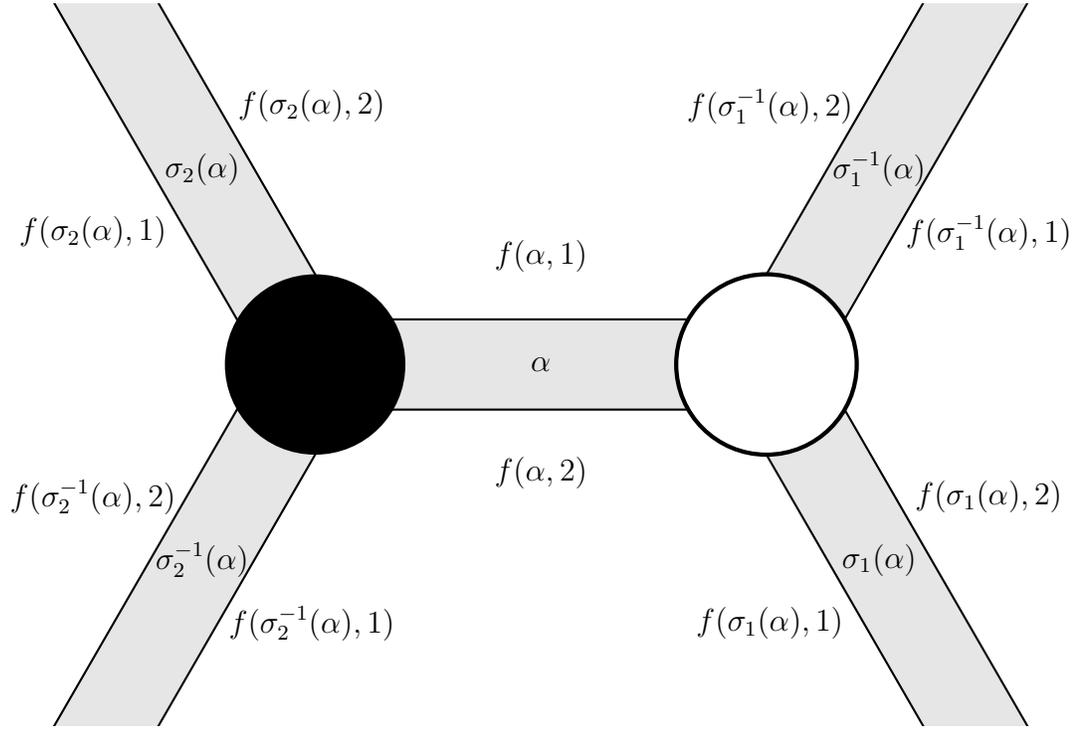
\begin{figure}[t]
\centering
\begin{tikzpicture}[scale=1.2]
\clip (-3.5,-4) rectangle (8.5,4);
\draw[draw=black,thick,fill=black!10] (0,-0.5) rectangle (5,0.5);
\begin{scope}[rotate=120]
   \draw[thick,fill=black!10] (0,-0.5cm) rectangle (10,0.5);
      \draw (2.5, 0) node {$\DELTA$};
      \draw (2.5, 1.4) node {$f(\DELTA,1)$};
      \draw (2.5,-1.4) node {$f(\DELTA,2)$};
\end{scope}
\begin{scope}[rotate=-120]
   \draw[thick,fill=black!10] (0,-0.5cm) rectangle (10,0.5);
      \draw (2.5, 0) node {$\EPSILON$};
      \draw (2.5, 1.4) node {$f(\EPSILON,1)$};
      \draw (2.5,-1.4) node {$f(\EPSILON,2)$};

\end{scope}

\begin{scope}[xshift=5cm]
    \begin{scope}[rotate=60]
      \draw[thick,fill=black!10] (0,-0.5) rectangle (10,0.5);
      \draw (2.5, 0) node {$\GAMMA$};
      \draw (2.5, 1.4) node {$f(\GAMMA,2)$};
      \draw (2.5,-1.4) node {$f(\GAMMA,1)$};
    \end{scope}
    \begin{scope}[rotate=-60]
      \draw[thick,fill=black!10] (0,-0.5) rectangle (10,0.5);
      \draw (2.5, 0) node {$\BETA$};
      \draw (2.5, 1.4) node {$f(\BETA,2)$};
      \draw (2.5,-1.4) node {$f(\BETA,1)$};

    \end{scope}
\end{scope}

\fill (0,0) circle (1); 
\draw[fill=white,ultra thick] (5,0) circle (1); 
\draw (2.5,1.2)  node{$f(\alpha,1)$};
\draw (2.5,-1.2) node {$f(\alpha,2)$};
\draw (2.5,0) node {$\alpha$};

\end{tikzpicture}
\caption{The convention for defining the labels of the edge-sides.
The symbol written directly on the edge is its label;
as usual the counterclockwise cycle structure of the edges around the white vertices
is given by the permutation $\sigma_1$ while the 
counterclockwise cycle structure of the edges around the black vertices
is given by the permutation $\sigma_2$.
The two symbols written next to the edge are the labels of the edge-sides.}
\label{fig:oriented-to-nonoriented}
\end{figure}

More formally speaking, the above construction corresponds to the 
non-oriented map $M=(\BC,\WC,\E)$, where
\begin{align*} 
\BC &= \Big\{ \big\{f(k,1), f\big( \sigma_2(k), 2 \big)\big\} : k \in [n] \Big\}, \\
\WC&=  \Big\{ \big\{f(k,2), f\big( \sigma_1(k), 1 \big)\big\} : k \in [n] \Big\}, \\  
\E &= \Big\{ \{f(k,1), f(k,2)\} : k \in [n] \Big\}. 
\end{align*}
The resulting non-oriented map $M$ is orientable and connected.
Furthermore, each such an orientable, connected, non-oriented map 
can be obtained in this way in $2 \cdot n!$ ways
(the factor $2$ comes from the two possible ways of choosing the orientation of the map,
the factor $n!$ counts the possible choices of the labeling of the edges of the oriented map
from which we start).

\medskip

The comparison of the multiplicities arising from the two ways of interpreting the summation
over \eqref{eq:transitive-pairs} concludes the proof.
\end{proof}

\section{Proof of the first main result, \cref{theo:main-bijection}}
\label{sec:proof:thm:main}

\begin{proof}[Proof of \cref{theo:main-bijection}]
Recall that the left-hand side of \eqref{eq:main-thm} runs over
oriented, unlabeled, rooted, connected maps $M_1$ with $n$ edges.
Each such a map can become a \emph{labeled} map in $(n-1)!$ distinct ways as follows:
we attach to the root edge the label $1$ and
we attach the remaining labels from $[n]\setminus\{1\}$ to the other edges. 
Thus the left-hand side of \eqref{eq:main-thm} is equal to
\begin{equation}
\label{eq:intermediate}
\frac{1}{(n-1)!} \sum_{\substack{\sigma_1,\sigma_2\in \Sym{n}, \\ 
\text{$\langle \sigma_1,\sigma_2\rangle$ is transitive}}} M(\sigma_1,\sigma_2)  
\end{equation}
which is viewed as linear combination of \emph{unlabeled} maps.

\medskip

Edge liberation for oriented maps (\cref{prop:edge-liberation-oriented}) 
implies that \eqref{eq:intermediate} is equal to
\begin{equation}
\label{eq:intermediate2}
 \frac{1}{(2n-1)!} \sum_M M, 
\end{equation}
where the sum runs over non-oriented maps
$M$ over the base set $X=[2n]$ such that $M$ is orientable and connected.
Our analysis of the left-hand side of \eqref{eq:main-thm} is now completed
and we shall turn to its right-hand side.

\bigskip

Recall that the sum on the right-hand side of \eqref{eq:main-thm}
runs over non-oriented, unlabeled, rooted maps $M_2$ with $n$ edges and one face.
Note that there is a canonical way of defining the labels on the edge-sides of this map,
as follows. We perform a walk along the boundary of the unique face of map,
with the first step along the decorated edge-side in the direction from its black
to its white extremity. We attach the label $1$ to this decorated edge-side,
and we continue to label the edge-sides in the order in which we visit them with
the remaining elements of the label set $[2n]\setminus\{1\}$.
This means that the sum on the right-hand side of \eqref{eq:main-thm}
is in fact a conservative summation over non-oriented maps with face-type $(n)$. 

\medskip

From \cref{def:mon-top-degree} it follows now that the right-hand side of
\eqref{eq:main-thm} is equal to 
\begin{equation}
\label{eq:inermediate-X1}
 \frac{1}{n!} \sum_{M=(G,S)} \sum_{\substack{\prec: \\ \text{$(M,\prec)$ is top-degree}}} G, 
\end{equation}
where the first sum is a \emph{conservative} sum over non-oriented maps with face-type $(n)$ 
and the second sum runs over the histories such that $(M,\prec)$ is a top-degree pair.
We apply edge liberation for non-oriented maps (\cref{prop:map-liberation-nonoriented}); 
in this way \eqref{eq:inermediate-X1} takes the form
\begin{equation}
\label{eq:inermediate-X2}
 \frac{1}{n!\ (2n-1)!} \sum_{M=(G,S)} \sum_{\substack{\prec: \\ \text{$(M,\prec)$ is top-degree}}} G,    
\end{equation}
where the first sum is this time a \emph{liberal} sum over non-oriented maps with face-type $(n)$.

We apply the bijection provided by \cref{theo:key-bijection}; in this way we see that \eqref{eq:inermediate-X2}
is equal to 
\begin{equation}
\label{eq:inermediate-X3}
 \frac{1}{(2n-1)!} \sum_{M=(G,S)}  G,    
\end{equation}
where the sum runs over all connected, orientable non-oriented maps over the base set $[2n]$.

\medskip

Since \eqref{eq:intermediate2} and \eqref{eq:inermediate-X3} are equal as formal linear combinations
of unlabeled graphs, this concludes the proof.
\end{proof}

\section{Background and history of the result: Jack characters}
\label{sec:introduction}

In this section we present the background and the motivations for \cref{theo:main-bijection},
in particular the problems related to \emph{Jack characters}.

\subsection{Jack characters}

\subsubsection{Jack polynomials}
\label{sec:jack-polynomials-motivations}
Henry Jack \cite{Jack1970/1971} introduced 
a family $\big( J^{(\alpha)}_\pi\big)$ (indexed by an integer partition $\pi$) of symmetric functions
which depend on an additional parameter $\alpha$. 
During the last forty years, 
many connections of these \emph{Jack polynomials} 
with various fields of mathematics and physics were established: 
it turned out that they play a crucial role in 
understanding Ewens random permutations model \cite{DiaconisHanlon1992}, 
generalized $\beta$-ensembles and some statistical mechanics models 
\cite{OkounkovOlshanski1997},
Selberg-type integrals \cite{Kaneko1993},
certain random partition models 
\cite{
BorodinOlshanski2005},
and some problems of algebraic geometry \cite{%
Okounkov2003}, 
among many others. 
Better understanding of Jack polynomials is also very desirable in the context of generalized 
$\beta$-ensembles and their discrete counterpart model \cite{DoleegaFeray2014}. 
Jack polynomials are a special case of the celebrated \emph{Macdonald polynomials} 
which \emph{``have found applications in special function theory, representation theory, 
algebraic geometry, group theory, statistics and quantum mechanics''} \cite{GarsiaRemmel2005}.

\subsubsection{Dual combinatorics of Jack polynomials}
Lassalle \cite{Lassalle2008a,Lassalle2009} initiated investigation of a kind of \emph{dual combinatorics} of 
Jack polynomials. 
More specifically, 
one expands Jack polynomial in the basis of power-sum symmetric functions:
\begin{equation} 
\label{eq:definition-theta}
J^{(\alpha)}_\lambda = \sum_\pi \theta^{(\alpha)}_\pi(\lambda)\ p_\pi. 
\end{equation}
The above sum runs over partitions $\pi$ such that $|\pi|=|\lambda|$. 
The coefficient $\theta^{(\alpha)}_\pi(\lambda)$ is called \emph{unnormalized Jack character};
with the right choice of the normalization it
becomes \emph{the normalized Jack character $\Ch_\pi(\lambda)$}
(the details of this relationship will be given in \cref{def:jack-character-classical}).
An interesting feature of Jack characters is that for the special choice
of the deformation parameter $\alpha=1$ they coincide with the usual characters of the symmetric groups.

The above approach is referred to as \emph{dual} because one fixes $\pi$ 
and views the character 
$\lambda\mapsto \Ch_\pi(\lambda)$ as a function of the Young diagram $\lambda$, 
opposite to the usual approach in the representation theory where one usually fixes 
the irreducible representation (the diagram $\lambda$) and views the character as a function of the conjugacy class 
(the partition $\pi$).
In the context of the representation theory of the symmetric groups 
(which corresponds to the special case $\alpha=1$)
this dual approach
was initiated by Kerov and Olshanski \cite{KerovOlshanski1994}
and it soon turned out to be highly successful;
its applications include, for example, random Young diagrams
\cite{IvanovOlshanski2002,Biane1998,'Sniady2006c}.
Lassalle \cite{Lassalle2008a,Lassalle2009} adapted this idea to the framework of Jack characters.

\smallskip

Jack characters $\Ch_\pi(\lambda)$ are in the focus of the current paper.
Our motivations for investigating them are threefold. 
Firstly, since Jack characters are related to Jack polynomials, a better understanding
of the former might shed some light on the latter. 
Secondly, they can be used in order to investigate some natural 
deformations of classical random Young diagrams
\cite{Sniady2016,DolegaSniady2014}.
Thirdly, numerical data \cite{LassalleData} as well as some partial theoretical results
\cite{Lassalle2008a,Lassalle2009,DolegaFeraySniady2013,Sniady2015,Sniady2016}
indicate that they might have a rich algebraic-combinatorial or
representation-theoretic structure.

\smallskip

Our ultimate goal would be to find some convenient
closed formula for Jack characters. This goal is beyond our reach;
a more modest goal would be to find a closed formula for
the dominant part of Jack characters in the suitable asymptotic scaling.
We shall address this issue later on.

\subsubsection{Definition of Jack characters}
In order for this dual approach to be successful (both with respect to the usual characters 
of the symmetric groups and for the Jack characters)
one has to choose the most convenient normalization constants.
In the current paper we will use the normalization introduced by Dołęga and F\'eray 
\cite{DoleegaFeray2014} 
which offers some advantages over the original normalization of Lassalle. 
Thus, with the right choice of the multiplicative constant, the unnormalized Jack character
$\theta_{\lambda}^{(\alpha)}(\pi)$ from \eqref{eq:definition-theta} 
becomes the \emph{normalized Jack character $\Ch^{(\alpha)}_\pi(\lambda)$}.
Regretfully, their definition is quite technical and not very enlightening.
On the bright side, this definition is not relevant for the purposes 
of the current paper and the Readers faint at heart are cordially invited to
fast forward to \cref{sec:definition-finished}.

\begin{definition}
\label{def:jack-character-classical}
Let $\alpha>0$ be given and let $\pi$ be a partition. The \emph{normalized Jack character
$\Ch_\pi(\lambda)$}  is given by:
\begin{equation}
\Ch_{\pi}(\lambda):=
\begin{cases}
\alpha^{-\frac{|\pi|-\ell(\pi)}{2}}
\binom{|\lambda|-|\pi|+m_1(\pi)}{m_1(\pi)}
\ z_\pi \ \theta^{(\alpha)}_{\pi,1^{|\lambda|-|\pi|}}(\lambda)
&\text{if }|\lambda| \ge |\pi| ;\\
0 & \text{if }|\lambda| < |\pi|,
\end{cases}
\end{equation}
where 
\[ z_\pi=\prod_i i^{m_i(\pi)}\ m_i(\pi)! \]
is the standard numerical factor.

Jack character $\Ch_{\pi}(\lambda)$ depends on the deformation parameter
$\alpha$, but to keep the notation light we shall usually make this dependence implicit.
\end{definition}

\subsubsection{The deformation parameter}
\label{sec:definition-finished}
In order to avoid dealing with the square root of the variable $\alpha$, 
we introduce an indeterminate $A$ such that
\[ A^2 = \alpha.\]
In this way Jack character $\Ch_{\pi}(\lambda)\in\Laurent$ 
becomes a Laurent polynomial in the indeterminate $A$.
This is the viewpoint which we will \emph{usually} take in this
paper, with the exception of \cref{sec:stanley1,sec:stanley212}
where $A$ will be some fixed real number.

\subsection{Bicolored graphs and their embeddings}

\begin{figure}
\centering
\subfloat[][]{\begin{tikzpicture}[scale=0.5,
white/.style={circle,draw=black,inner sep=4pt},
black/.style={circle,draw=black,fill=black,inner sep=4pt},
connection/.style={draw=black!80,black!80,auto}
]
\footnotesize

\begin{scope}
\clip (0,0) rectangle (10,10);

\draw (3.333,2.333) node (b1)    [black,label=110:$\textcolor{\kolorPi}{\Pi}$] {};
\draw (b1) +(10,0) node (b1prim) [black] {};

\draw (7.666,6.666) node (b2)     [black,label=0:\textcolor{\kolorSigma}{$\Sigma$}] {};
\draw (b2) +(0,-10) node (b2prim) [black] {};

\draw (b2) +(-3,1) node (w2) [white,label=180:$\textcolor{\kolorW}{W}$] {};

\draw (6.666,3.333) node (w1) [white,label=120:$\textcolor{\kolorV}{V}$] {};
\draw (w1) +(-10,0) node (w1left) [white] {};
\draw (w1) +(0,10)  node (w1top)  [white] {};

\draw[connection] (b1) to node {\textcolor{black}{$4$}} node [swap] {} (w1);

\draw[connection] (b2) to node {\textcolor{black}{$3$}} node [swap] {} (w2);

\draw[connection,pos=0.333] (b2) to node {\textcolor{black}{$2$}} node [swap] {} (w1top);
\draw[connection,pos=0.666] (b2prim) to node {\textcolor{black}{$2$}} node [swap] {} (w1);

\draw[connection,pos=0.666] (b1prim) to node {\textcolor{black}{$1$}} node [swap] {} (w1);
\draw[connection,pos=0.333] (b1) to node {\textcolor{black}{$1$}} node [swap] {} (w1left);

\draw[connection,pos=0.666] (w1) to node {} node [swap] {\textcolor{black}{$5$}} (b2);

\end{scope}

\draw[very thick,decoration={
    markings,
    mark=at position 0.5 with {\arrow{>}}},
    postaction={decorate}]  
(0,0) -- (10,0);

\draw[very thick,decoration={
    markings,
    mark=at position 0.5 with {\arrow{>}}},
    postaction={decorate}]  
(0,10) -- (10,10)  ;

\draw[very thick,decoration={
    markings,
    mark=at position 0.5 with {\arrow{>>}}},
    postaction={decorate}]  
(0,0) -- (0,10);

\draw[very thick,decoration={
    markings,
    mark=at position 0.5 with {\arrow{>>}}},
    postaction={decorate}]  
(10,0) -- (10,10)  ;
\end{tikzpicture}
\label{subfig:map}}
\hspace{5ex}
\subfloat[][]{
\begin{tikzpicture}[scale=1.2]
\begin{scope}

\draw[line width=5pt,\kolorSigma!20] (-0.2,0.5) -- (3.2,0.5);
\draw (3.2,0.5) node[anchor=west] {$\textcolor{\kolorSigma}{\Sigma}$};

\draw[line width=5pt,\kolorPi!20] (-0.2,1.5) -- (1.2,1.5);
\draw (1.2,1.5) node[anchor=west] {$\textcolor{\kolorPi}{\Pi}$};

\draw[line width=5pt,\kolorW!20] (2.5,-0.2) -- (2.5,1.2);
\draw (2.5,1.2) node[anchor=south] {$\textcolor{\kolorW}{W}$};

\draw[line width=5pt,\kolorV!20] (0.5,-0.2) -- (0.5,2.2);
\draw (0.5,2.2) node[anchor=south] {$\textcolor{\kolorV}{V}$};

\draw[ultra thick] (4.5,0) -- (3,0) -- (3,1) -- (1,1) -- (1,2) -- (0,2) -- (0,2.5); 
\draw (0,0) -- (3,0) -- (3,1) -- (1,1) -- (1,2) -- (0,2); 
\clip (0,0) -- (3,0) -- (3,1) -- (1,1) -- (1,2) -- (0,2); 
\draw (0,0) grid (3,3);
\end{scope}
\draw (0.5,-0.2) node[anchor=north,text height=8pt] {$a$};
\draw (1.5,-0.2) node[anchor=north,text height=8pt] {$b$};
\draw (2.5,-0.2) node[anchor=north,text height=8pt] {$c$};

\draw (-0.2,0.5) node[anchor=east]  {$\alpha$};
\draw (-0.2,1.5) node[anchor=east] {$\beta$};

\draw(2.5,0.5) node {$3$};
\draw(0.5,0.5) node {$2,5$};
\draw(0.5,1.5) node {$1,4$};

\end{tikzpicture}
\label{subfig:embed}}

\caption{\protect\subref{subfig:map} Example of a bicolored graph (drawn on the torus: the left side of the square 
should be glued to the right side, as well as bottom to top, as indicated by arrows) and 
\protect\subref{subfig:embed}  an example of its embedding 
$F(\Sigma)=\alpha$, $F(\Pi)=\beta$, $F(V)=a$, $F(W)=c$.
$F(1)=F(4)=(a \beta)$, $F(2)=F(5)=(a \alpha)$, $F(3)=(c\alpha)$.
The columns of the Young diagram were indexed by small Latin letters, the rows by small Greek letters.}
\label{fig:embedding}
\end{figure}
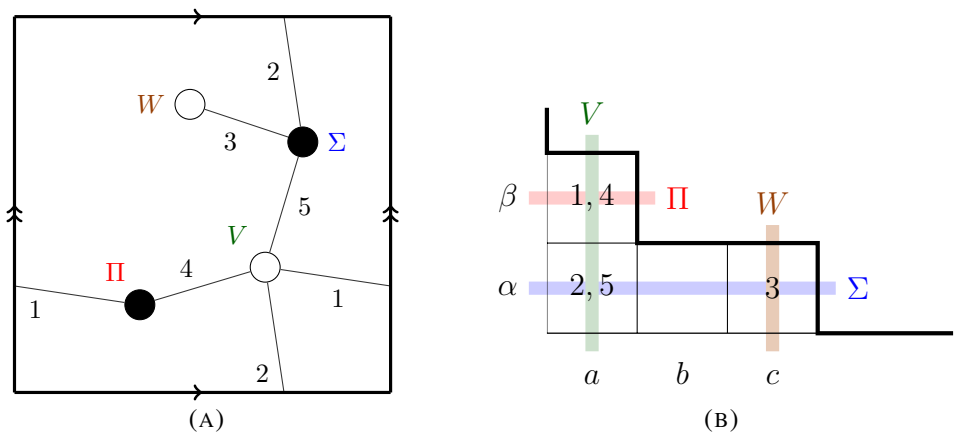

An \emph{embedding $F$ of a bicolored graph $G$ to a Young diagram $\lambda$} is a function which maps the set 
$V_\circ(G)$ of white vertices of $G$ to the set of columns of $\lambda$, which maps the set $V_\bullet(G)$ of 
black vertices of $G$ to the set of rows of $\lambda$, and maps the edges of $G$ to boxes of $\lambda$, see \cref{fig:embedding}. We also require that an embedding preserves the relation of \emph{incidence},
i.e.,~a vertex $V$ and an incident edge $E$ should be mapped to a row or column $F(V)$ which contains the box
$F(E)$. We denote by $N_G(\lambda)$ the number of such embeddings of $G$ to $\lambda$. 
Quantities $N_G(\lambda)$ were introduced by F\'eray and the second-named author \cite{FeraySniady2011a} 
and they proved to be very useful for studying various asymptotic and enumerative problems of the 
representation theory of symmetric groups \cite{FeraySniady2011a,FeraySniady2011,DolegaF'eray'Sniady2008}.

\begin{definition}
For a bicolored graph $G$ and a Young diagram $\lambda$
we define the \emph{normalized number of embeddings} \cite{Sniady2015} which is a Laurent polynomial in $A$:
\begin{equation}
\label{eq:normalized-embedding}
 \Embed_G (\lambda):=
\frac{{A}^{|V_\circ(G)|}}{(-A)^{|V_\bullet(G)|}} 
\ N_{G}(\lambda) \in\Laurent.
\end{equation}
\end{definition}
A very similar quantity denoted by $N_G^{(\alpha)}(\lambda)$ --- 
which differs from $\Embed_G (\lambda)$ only by the choice of the sign --- 
was considered already in \cite{DolegaFeraySniady2013}.

\medskip

In the case when $M=(G,S)$ is a map, we denote
$\Embed_M (\lambda):= \Embed_G (\lambda)$.

\subsection{Stanley formulas}
\label{sec:stanley}

For some special values of the deformation parameter 
$A\in\left\{\pm \frac{1}{\sqrt{2}}, \pm 1, \pm \sqrt{2} \right\}$
which correspond to
$\alpha\in\left\{\frac{1}{2},1,2\right\}$,
Jack characters admit closed formulas in terms of 
\emph{embeddings} of certain \emph{bicolored maps}.
Formulas of this type are called \emph{Stanley formulas} after
Richard Stanley, who found such a formula for $\alpha=1$
as a conjecture~\cite{Stanley-preprint2006}.
We shall review them in the following.

\subsubsection{Stanley formula for $\alpha=1$ and oriented maps}
\label{sec:stanley1}

In the special case of $A= \pm 1$ which corresponds to $\alpha=1$, 
Jack polynomials coincide
(up to simple multiplicative constants) with Schur polynomials.
Using this fact one can show that 
in this special case the Jack character $\Ch_\pi^{A:= 1}$ coincides 
with the (suitably normalized) character of the symmetric group;
for the details see the work of Lassalle \cite{Lassalle2009} (who used a different normalization)
as well as the work of Dołęga and F\'eray \cite{DoleegaFeray2014}.
For this reason, for $A=1$ Jack characters have a much richer algebraic and 
representation-theoretic structure than for a generic value of $A$.

\medskip

In particular, it has been observed in \cite{FeraySniady2011a} that a certain formula 
conjectured by Stanley \cite{Stanley-preprint2006}
and proved by F\'eray \cite{F'eray2010}
for the normalized characters of the symmetric groups 
can be expressed as the sum
\begin{equation}
\label{eq:alpha=1}
 \Ch^{A:=1}_\pi(\lambda) =(-1)^{\ell(\pi)} \sum_{M}  \Embed^{A:=1}_M(\lambda) 
\end{equation}
over all \emph{oriented bicolored maps $M$ with face-type $\pi$},
where $\ell(\pi)$ denotes the number of parts of the partition $\pi$.

\subsubsection{Stanley formula for $\alpha\in\{2,\frac{1}{2}\}$ and non-oriented maps}
\label{sec:stanley212}
In the special case when $A=\pm \sqrt{2}$ and $\alpha=2$ (respectively, 
$A=\pm \frac{1}{\sqrt{2}}$ and $\alpha=\frac{1}{2}$) 
Jack polynomials coincide with \emph{zonal polynomials}
(respectively, \emph{symplectic zonal polynomials}).
Thanks to this additional structure it has been proved in 
a joint work of the second-named author with F\'eray \cite{FeraySniady2011} that
\begin{align} 
\label{eq:alpha=2}
\Ch_\pi^{{A:=\sqrt{2}}} (\lambda) &= (-1)^{\ell(\pi)}\ \sum_M  
\left(- \frac{1}{\sqrt{2}} \right)^{|\pi|+\ell(\pi)-|\V(M)|}
  \Embed^{A:=\sqrt{2}}_M(\lambda),\\
\label{eq:alpha=1/2}
\Ch_\pi^{{A:=\frac{1}{\sqrt{2}}}} (\lambda) &= (-1)^{\ell(\pi)}\ \sum_M  
\left(\frac{1}{\sqrt{2}} \right)^{|\pi|+\ell(\pi)-|\V(M)|}
  \Embed^{A:=1/\sqrt{2}}_M(\lambda), 
\end{align}
where the sums run over all \emph{non-oriented maps $M$
with face-type $\pi$}, as in \cref{sec:summation-face-type}.
The reader should be advised that
the notations and the normalizations in \cite{FeraySniady2011}
are a bit different;
the link between the statements above and the results of \cite{FeraySniady2011}
is given in \cite[Section 5]{DolegaFeraySniady2013}.

\subsection{How to prove a closed formula for Jack characters?
Orientability generating series}
\label{subsec:how-to-prove}

Jack characters admit an abstract characterization \cite[Theorem 1.7]{Sniady2015}
which was found in a recent paper of the second-named author
(see also Theorem A.2 of F\'eray in an appendix to the same paper \cite{Sniady2015}).
This abstract characterization opens the following path toward \emph{proving}
a closed formula for Jack characters: in the first step one should \emph{guess}
the right formula, and in the second step one should \emph{verify} that it indeed 
fulfills the aforementioned defining properties of Jack characters.

\medskip

How to make the first step and to guess a closed formula for Jack characters?
The three Stanley formulas \eqref{eq:alpha=1}, \eqref{eq:alpha=2},  \eqref{eq:alpha=1/2}
might suggest that the hypothetical formula for Jack characters in the generic case should be of the form
\begin{equation}
\label{eq:wow-conjecture}
 \Ch_\pi(\lambda) = (-1)^{\ell(\pi)} \sum_M \operatorname{weight}_M  \Embed_M(\lambda),    
\end{equation}
where the sum should run over \emph{non-oriented} maps with face-type $\pi$.
In the above formula $\operatorname{weight}_M\in\Laurent$ is some hypothetical quantity which
measures the \emph{non-orientability} of the map $M$ 
\cite[Conjecture 1.1]{DolegaFeraySniady2013}.

\medskip

A joint work of the second-named author with Dołęga and F\'eray \cite{DolegaFeraySniady2013}
presents an attempt to guess the exact form of this hypothetical quantity $\operatorname{weight}_M$,
an attempt which was based on a reverse-engineering of the results of Lassalle \cite{Lassalle2008a}.
Our candidate quantity, denoted $\weight_M$ (which is an acronym for 
\emph{\textbf{m}easure \textbf{o}f \textbf{n}on-orientability})
was defined as follows:
we remove the edges from the ribbon graph of $M$ in a uniformly random order;
to each edge which is about to be removed we associate a factor which is
related to the topological way in which this edge is attached to the remaining 
edges (i.e., the edges which have not been removed yet) as we discussed in
\cref{subsec:anatomy}. 
The quantity $\weight_M$ is defined
as the expected value of the product of the aforementioned factors.
The details of this construction will be recalled in \cref{subsec:measure-of-nonorientability}.
This weight $\weight_M$ gives rise to the \emph{orientability generating series}
\cite[Section 1.10]{DolegaFeraySniady2013} which is defined in analogy to \eqref{eq:wow-conjecture} as
\begin{equation}
\label{eq:what-is-ogs}
 \widehat{\Ch}_\pi(\lambda) := (-1)^{\ell(\pi)} \sum_M \weight_M  \Embed_M(\lambda),
\end{equation}
where the sum runs over \emph{non-oriented} maps with face-type $\pi$.

\subsection{Stanley polynomials}
In order to be able to speak about the \emph{degree}
of some functions on the set of Young diagrams we will need 
the notion of \emph{Stanley polynomials}.
The content of this section is an abridged and less formal version of 
\cite[Section 1.10]{Sniady2015}.

\begin{figure}[t]
\begin{tikzpicture}[scale=0.8]

\begin{scope}
\clip (0,0) -- (11,0) -- (11,2) -- (8,2) -- (8,5) -- (4,5) -- (4,7) -- (0,7) -- cycle;
\draw[black!60] (0,0) grid[step=0.5] (20,10); 
\end{scope}

\draw[ultra thick](0,0) -- (11,0) -- (11,2) -- (8,2) -- (8,5) -- (4,5) -- (4,7) -- (0,7) -- cycle;

\draw[dotted,thick] (8,2) -- (0,2)
(4,5) -- (0,5);

\begin{scope}[<->,thick,auto,dashed,blue]
\draw (11.3,0) to node[swap] {$p'_1$} (11.3,2);
\draw (0,1.7) to node[swap,shape=circle,fill=white] {$q'_1$} (11,1.7);

\draw (8.3,2) to node[swap] {$p'_2$} (8.3,5);
\draw (0,4.7) to node[swap,shape=circle,fill=white] {$q'_2$} (8,4.7);

\draw (4.3,5) to node[swap] {$p'_3$} (4.3,7);
\draw (0,6.7) to node[swap,shape=circle,fill=white] {$q'_3$} (4,6.7);
\end{scope}

\end{tikzpicture}

\caption{Multirectangular Young diagram $P\times Q$.}

\label{fig:multirectangular}
\end{figure}
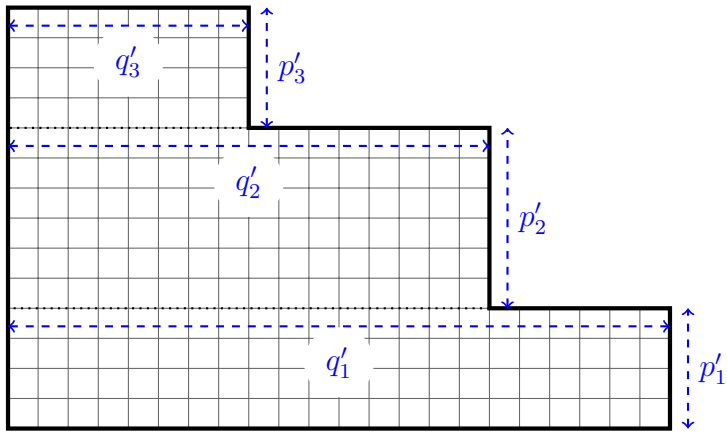

\subsubsection{Multirectangular coordinates}

We start with \emph{anisotropic multirectangular coordinates} $P=(p_1,\dots,p_\ell)$ and $Q=(q_1,\dots,q_\ell)$.
They give rise to \emph{isotropic multirectangular coordinates} given by
\begin{align*}
P'=(p'_1,\dots,p'_\ell):&= \left( A p_1, \dots, A p_\ell \right),\\
Q'=(q'_1,\dots,q'_\ell):&= \left( \frac{1}{A} q_1, \dots, \frac{1}{A} q_\ell \right).
\end{align*}
Note that $P'$ and $Q'$ depend implicitly on $P$ and $Q$.

Suppose that $P'=(p'_1,\dots,p'_\ell)$ and $Q=(q'_1,\dots,q'_\ell)$ are sequences of non-negative integers 
such that $q'_1\geq \cdots \geq q'_\ell$; we consider the \emph{multirectangular} Young diagram
\[P' \times Q' = (\underbrace{q'_1,\dots,q'_1}_{\text{$p'_1$ times}},\dots,
\underbrace{q'_\ell,\dots,q'_\ell}_{\text{$p'_\ell$ times}}).\]
This concept is illustrated in \cref{fig:multirectangular}.

\subsubsection{The deformation parameter $\gamma$}
We usually view Jack polynomials as functions of the parameter $\alpha$
and Jack characters as functions of another  parameter $A=\sqrt{\alpha}$.
However, it is convenient to consider yet another deformation parameter
\begin{equation}
\label{eq:gamma}
 \gamma:=\frac{1}{A}-A \in \Laurent.
\end{equation}
As we shall see, several quantities can be expressed as polynomials in $\gamma$.

\subsubsection{Stanley polynomials and degree of functions on $\AllYoung$}

Let $F\colon \AllYoung \to\Laurent$ be a function on the set of Young diagrams and 
let
\[\St_\ell=\St_\ell(\gamma;p_1,\dots,p_\ell;q_1,\dots,q_\ell)=\St_\ell(\gamma; P; Q) \] 
be a polynomial in $2\ell+1$ variables.

Suppose that the equality
\[  F(P'\times Q')=\St_\ell\left(\gamma ;P ; Q\right)\]
--- with the substitution \eqref{eq:gamma}
for the variable $\gamma$ --- 
holds true for all choices of $P$, $Q$ and $A\neq 0$ for which the multirectangular diagram $P'\times Q'$ is well-defined.
Then we say that $\St_\ell$ is the \emph{Stanley polynomial} for $F$.
The above definition is not very precise; 
on the formal level one should consider Stanley polynomial as an element of some inverse limit for $\ell\to\infty$;
for the details we refer to \cite[Section 1.10]{Sniady2015}.

Stanley polynomials are a perfect tool for studying asymptotic questions
in the setup when the Young diagram $P'\times Q'$ tends to infinity.
In particular, we say that $F\colon \AllYoung \to\Laurent$ 
is a function of degree (at most) $d$ if the corresponding
Stanley polynomial is of degree (at most) $d$.

\subsubsection{Stanley polynomials for $\Ch_\pi$ and $\widehat{\Ch}_\pi$}
It is a highly non-trivial result of Dołęga and F\'eray 
\cite[Corollary 3.5]{DoleegaFeray2014} 
(see also \cite[Theorem~2.15, Corollary~2.11]{Sniady2015} for details how to
adapt their result to our setup)
that a Stanley polynomial exists for each Jack character $\Ch_\pi$
and that is of degree at most $|\pi|+\ell(\pi)$. 
For example,
\begin{align*}
\Ch_1(P'\times Q') & = \sum_i  p_i q_i, 
\notag \\
\Ch_2(P'\times Q')  & = \sum_i 
 p_i q_i \left[  q_i - p_i+\gamma \right]
- 2 \sum_{i<j}  p_i p_j q_j, 
\notag \\
\Ch_3(P'\times Q') & =
\sum_i 
p_i q_i  \left[
q_i^2 -
3 p_i q_i 
+  p_i^2 +
3 \gamma (q_i- p_i) 
  +
 2 \gamma^2  + 1 \right] + \\
\notag
& -3 \sum_{i<j}  p_i p_j q_j
\left[ (q_i-p_i+\gamma)+(q_j-p_j+\gamma)  \right] + \\  
\notag
 &+ \sum_{i<j<k} 6 p_i p_j p_k q_k.
\end{align*}

One can show that the Stanley polynomial for 
the orientability generating series $\widehat{\Ch}_\pi$ exists
and it is also of degree $|\pi|+\ell(\pi)$;
the proof is postponed to \cref{lem:degree-bound}.

\subsection{Orientability generating series versus Jack character}
\label{sec:heuristic}
The initial prediction of the authors of \cite{DolegaFeraySniady2013} was that
the quantity $\weight_M$ is the right guess for the hypothetical $\operatorname{weight}_M$
and thus $\widehat{\Ch}_\pi = \Ch_\pi$.
Regretfully, this turns out to be not the case \cite[Section 7]{DolegaFeraySniady2013}
and $\widehat{\Ch}_\pi \neq \Ch_\pi$ in general.

This might seem as the end of the story and an example of a failed research,
nevertheless the orientability series $\widehat{\Ch}_\pi$ appeared to predict the 
properties of the Jack character $\Ch_\pi$ \emph{suspiciously well}.
The latter statement was supported heuristically in the following two ways.
\begin{itemize}
   \item 
Firstly, 
computer-assisted comparison of the coefficients of Stanley polynomials
for the character $\Ch_\pi$ with their counterpart for the orientability generating series
$\widehat{\Ch}_\pi$
indicates that \emph{a lot of them} coincide for several concrete examples.

\item Secondly, in the case when $\pi=(n)$ is a partition with only one part, 
a comparison of
\begin{enumerate}[label=(\emph{\roman*})]
   \item \label{item:A}
the contribution to the orientability generating series $\widehat{\Ch}_n$ 
which comes from the maps $M$ with small genus $g\leq \frac{3}{2}$ with 
   \item \label{item:B} its counterpart for Jack character $\Ch_n$
\end{enumerate}
indicate a full match \cite[Section 5]{ACJ2016}. 
For \ref{item:A} a closed formula found by the first named author \cite{ACJ2016} is available.
The comparison can be either performed numerically for small values of $n$
using the data provided by Lassalle \cite{LassalleData}, or
by comparing the aforementioned closed formula for \ref{item:A}
with Lassalle's conjectural closed formula  \cite[Section 11]{Lassalle2009};
they turn out to be identical.
\end{itemize}

\subsection{The second main result: the top-degree of Jack character}

One of the main results of the current paper is the following positive result.

\begin{theorem}[The second main result]
\label{theo:top-degre-equal}
For each integer $n\geq 1$
\begin{equation}
\label{eq:chtops-are-born-equal}
 \Ch_n^{\ttop} = \widehat{\Ch}_n^{\ttop}. 
\end{equation}
\end{theorem}

Above, $\Ch_n^{\ttop}$ denotes the \emph{top-degree of Jack character $\Ch_n$}.
More explicitly, Stanley polynomial for $\Ch_n$ is known to be of degree $n+1$;
the homogeneous part of this Stanley polynomial of degree $n+1$
defines a function on the set of Young diagrams which will be denoted $\Ch_n^{\ttop}$.
Analogously, the top-degree $\widehat{\Ch}_n^{\ttop}$ of the orientability generating series 
is defined as the homogeneous part of $\widehat{\Ch}_n$ of degree $n+1$.

\medskip

In the remaining part of this section we will explain the relationship
between the two main results of the current paper:
\cref{theo:main-bijection} and \cref{theo:top-degre-equal}.

\subsection{Two formulas for the top-degree part of Jack characters}

The following result was proved by the second-named author using 
a modified version of the strategy which we outlined in \cref{subsec:how-to-prove}:
guess the right closed formula and then verify that it satisfies
some abstract characterization (this abstract characterization
of $\Chtt_n$ turns out to be much more complex than the analogous
characterization of $\Ch_n$).

\begin{theorem}[{\cite[Theorem 1.21]{Sniady2015}}]
\label{theo:formula-for-chtop}
For each integer $n\geq 1$
\begin{equation}
\label{eq:formula-for-chtop}
{\Ch}_n^{\ttop} (\lambda) = (-1) 
\sum_M \gamma^{n+1-|\V(M)|}\ \Embed_M(\lambda),
\end{equation}
where the sum runs over 
\emph{oriented, unlabeled, rooted, connected maps with $n$ edges and with arbitrary face-type}.
\end{theorem}

Our proof of \cref{theo:top-degre-equal} will be based on the following
simple observation: our main bijective result (\cref{{theo:main-bijection}})
shows equality between the right-hand sides of 
\eqref{eq:what-is-ogs} and \eqref{eq:formula-for-chtop}.
The missing details of the proof will be provided in \cref{sec:proof-second-main-result}.

\subsection{Top-degree of Jack characters: the real story.
Motivations for \cref{theo:main-bijection}}
\label{sec:motivations1}

The (yet unpublished at the time) bijective result of the current paper (\cref{theo:main-bijection})
preceded the proof of the above-mentioned closed formula for the top-degree of Jack characters
from \cref{theo:formula-for-chtop}.
Indeed, the key difficulty in \cref{theo:formula-for-chtop} was \emph{to guess the right formula}
and this difficulty was overcome by
converting a conjectural formula based on the ideas from \cite{DolegaFeraySniady2013}
to a more convenient form thanks to \cref{theo:main-bijection}.
For more discussion on this topic see below.

\subsection{Outlook into the future}
\label{sec:motivations2}
The Reader interested in Jack characters may wonder: 
why bother now proving the second main result of the current paper,
\cref{theo:top-degre-equal}.
After all, a simpler closed formula for $\Chtt_n$ (\cref{theo:formula-for-chtop}) 
was already available;
the formula \eqref{eq:chtops-are-born-equal} has served its duty as
a source of heuristics and can be retired now.
 
However, our ultimate goal, finding a closed formula for Jack characters $\Ch_n$
(and, more generally, $\Ch_\pi$), has not been completed. 
The latter is currently far beyond our reach;
a more modest partial accomplishment would be to find a closed formula for
the next term in the asymptotic expansion of Jack character $\Ch_n$ after its top-degree
part $\Ch^{\ttop}_n$. Once the right candidate formula is found by some heuristic means,
the machinery from \cite{Sniady2015} could  be probably relatively easily adapted 
in order to prove that such a candidate formula indeed holds true. 

Regretfully, there is a missing key element to this approach: we have no good candidate
formula for such a sub-dominant part of $\Ch_n$ of degree $n-1$
(see \cite[Section 1.17]{Sniady2015}). 
Probably the simplest approach would be to deduce or extrapolate such a formula based
on existing formulas for the asymptotically dominant part $\Ch_n^{\ttop}$.
Unfortunately, the final formula \eqref{eq:formula-for-chtop} for $\Ch_n^{\ttop}$
does not seem to offer any hints how to extrapolate it into the sub-dominant regime.

On the other hand, our original starting point, the orientability generating 
series $\widehat{\Ch}_n$ given by \eqref{eq:what-is-ogs}, offers immediately
a candidate formula for such a sub-dominant part. We have to admit:
it is a candidate formula which gives \emph{slightly} wrong predictions
\cite[Section 7]{DolegaFeraySniady2013}, but it is nevertheless a good 
starting point for some better formula.

Suppose that this happens to be indeed the case and someone, someday finds such a
better formula based on the concept of the orientability generating series.
In such a scenario some modified version of the bijection behind our main bijective result
(\cref{theo:main-bijection})
might become handy in order to convert this hypothetical formula into a more
convenient form, just like we used this bijection in the current paper in order 
to transform 
the right-hand side of \eqref{eq:what-is-ogs} 
into more convenient the right-hand side of \eqref{eq:formula-for-chtop}.

\medskip

For more discussion on this topic see \cite[Section 1.17]{Sniady2015}.

\subsection{$b$-conjecture}

In this section all maps are \emph{connected}, and \emph{rooted}, 
that is they posses a marked, oriented corner incident to some black vertex 
(that is, an angular region around black vertex $v$ delimited by two consecutive edges 
attached to $v$).

Goulden and Jackson \cite{GouldenJackson1997} introduced, using Jack symmetric functions, some multivariate generating series $\psi(\bm{x}, \bm{y}, \bm{z}; 1, 1+\beta)$ with an additional parameter $\beta$ that might be 
interpreted as a continuous deformation of the rooted bicolored maps generating series.
Indeed, it has the property that for $\beta \in \{0,1\}$
it specializes to rooted, orientable (for $\beta=0$) or
general, i.e.~orientable or not (for $\beta=1$) bicolored maps generating series.
Goulden and Jackson made the following conjecture:
coefficients of $\psi$ (indexed by three partitions of the same size) are polynomials in $\beta$ with positive integer coefficients that can be written as a multivariate generating series of rooted, general bicolored maps,
where the exponent of $\beta$ is an integer-valued statistics
that in some sense ``measures the non-orientability'' of
the corresponding bicolored map.

This \emph{$b$-conjecture} has not yet been proved, 
but some progress towards determining both algebraic 
and combinatorial properties of the coefficients in question has been made, 
and the work is ongoing. 
Dołęga and F\'eray \cite{DolegaFeray2016} proved recently that all coefficients of $\psi$ 
are polynomials in $\beta$ with rational coefficients. 
Based on this result, Dołęga \cite{Dolega2016} 
found a combinatorial interpretation of the top-degree of coefficients of $\psi$ 
indexed by two arbitrary partitions $\mu,\nu \vdash n$, 
and one partition consisting of only one part $(n)$, 
which conjecturally should be given by certain maps with only one face. 
An interesting phenomenon is that, similarly as in our \cref{theo:main-bijection}, 
the top-degree part found by Dołęga is given by 
\emph{orientable} maps with the black, and white, 
respectively, vertex degrees given by partitions $\mu$, and $\nu$, 
respectively, and \emph{arbitrary face structure}, and at the same time 
it is given by certain maps (called \emph{unhandled}) 
with a unique face, and black and white 
vertex degrees given by partitions $\mu$, and $\nu$, respectively.

We cannot resist to state that there must be a deep connection between these problems, and understanding it would be of great interest.

\section{Top-degree of the orientability generating series}
\label{section:clever-appendix}
\label{sec:final-section}

In this section we will present the details of the definition of the polynomial
$\weight_M(\gamma)$, its relationship to the quantity $\weight^{\ttop}_{M}$
introduced in \cref{def:mon-top-degree},
and the missing details of the proof of \cref{theo:top-degre-equal}.

\subsection{Weight associated to a map with a history}
\label{sec:weight-map-history}
We continue the discussion from \cref{sec:histories}.

Let $E_1,\dots,E_n$ be the sequence of edges of a non-oriented $M$, 
listed according to the linear order~$\prec$. We set 
$M_i = M \setminus \{ E_1,\dots, E_i \}$ and define
\begin{equation}
\label{eq:weight-map-order}
 \weight_{M,\prec} := \prod_{0\leq i \leq n-1}  \weight_{M_{i}, E_{i+1}}.  
\end{equation}
This quantity \eqref{eq:weight-map-order} 
can be interpreted as follows: from the map $M$ we remove (one by
one) all the edges, in the order
specified by the history.
For each edge which is about to be removed we consider
its weight $\weight_{M_{i}, E_{i+1}}$ 
relative to the current map
(recall that the factor $\weight_{M_{i}, E_{i+1}}$
was defined in \cref{subsec:anatomy}
and it depends on the \emph{type} of the edge $E_{i+1}$ in the map $M$,
i.e.~\emph{straight} versus \emph{twisted} versus \emph{interface}).

\subsection{The top-degree of $\weight_{M,\prec}$}
\label{sec:top-degree-of-weight}

The following result provides some crude information about the 
polynomial $\weight_{M,\prec}(\gamma)$.
\begin{lemma}[{\cite[Lemma 3.7]{DolegaFeraySniady2013}}]
\label{lem:how-twisted}
For any map $M$ the weight 
$\weight_{M,\prec}$ is a polynomial in the variable $\gamma$ of degree (at most)
$ 2\genus(M)$. 
\end{lemma}

Here we use the term \emph{genus} with a small abuse of notation; 
usually it is used only for orientable surfaces while we use it also for a 
\emph{non-orientable} connected map by setting 
\[\genus(M):=\frac{2 \cdot (\text{the number of connected components of $M$})-\chi(M)}{2},\]
where 
\[ \chi(M)= |\F(M)| - |\E(M)| + |\V(M)| \]
is the Euler characteristic of $M$. 

\medskip

It follows that the degree of the polynomial $\weight_{M,\prec}$
is bounded from above by 
\begin{multline}
\label{eq:DFSbound-prim}
\operatorname{degree} \weight_{M,\prec} \leq \\ 
 2 \ |\F(M)| -\chi(M) = |\F(M)|+ |\E(M)| - |\V(M)|.   
\end{multline}

The remaining part of this section is devoted to the investigation of the 
corresponding leading coefficient
\[ \left[ \gamma^{|\F(M)|+ |\E(M)| - |\V(M)|} \right] \weight_{M,\prec}(\gamma).\]

\bigskip

The following result is an extension of \cref{lem:connected-component-one-face-abbridged};
its only new component is condition \ref{enum:topdegree}.

\begin{lemma}
\label{lem:connected-component-one-face}
Let $M$ be a non-oriented map and $\prec$ be a history.
We use the notations from \cref{sec:histories}, i.e.~$E_1,\dots,E_n$ are the sequence of edges of $M$, 
listed according to the
linear order~$\prec$ and $M_i = M \setminus \{ E_1,\dots, E_i \}$.

Then the following conditions are equivalent:
\begin{enumerate}[label={(\Alph*)}]
   \item \label{enum:maptopdegree}
the pair $(M,\prec)$ is top-degree (\cref{def:top-degree});

   \item \label{enum:types-of-edges}
      for each $0\leq i\leq n-1$ the edge $E_{i+1}$ of the map $M_i$ is either:
     \begin{itemize}
        \item a twisted edge, or,
        \item a bridge or a leaf;
     \end{itemize}

   \item the degree of the polynomial $\weight_{M,\prec}$ is \emph{equal} to $|\F(M)|+ |\E(M)| - |\V(M)|$.
         \label{enum:topdegree} 
\end{enumerate}

\medskip

If the above conditions hold true then the leading coefficient of the polynomial
$\weight_{M,\prec}$
is given by
\begin{equation}
\label{eq:monic}
 \left[ \gamma^{|\F(M)|+ |\E(M)| - |\V(M)|}\right] \weight_{M,\prec} = 1.    
\end{equation}
\end{lemma}
\begin{proof}
Condition \ref{enum:topdegree} holds true if and only if
the inequalities involved in the proof of the bound \eqref{eq:DFSbound-prim}
become equalities and this happens
when both of the following two conditions hold true:
\begin{enumerate}[label={(C\arabic*)}]
   \item \label{cond:oneface} the map $M$ is top-degree (\cref{def:top-degree-map}), and,
   \item \label{cond:weight-max}
the weight $\weight_{M,\prec}$ is a polynomial in the variable $\gamma$ of degree \emph{exactly} $2\genus(M)$.
\end{enumerate}
In the following we shall find some equivalent reformulations of the condition \ref{cond:weight-max}.

\medskip

By revisiting the proof of \cref{lem:how-twisted} presented in \cite[Lemma 3.7]{DolegaFeraySniady2013}
one can show that the condition \ref{cond:weight-max} holds true
if and only if for each $0\leq i\leq {n-1}$ 
the weight $\weight_{M_i,E_{i+1}}$ is a polynomial in the variable $\gamma$ of degree \emph{exactly} 
\[2\genus(M_i)-2\genus(M_i\setminus E_{i+1})=2\genus(M_i)-2\genus(M_{i+1}).\]
The case-by-case analysis from the proof of \cite[Lemma 3.7]{DolegaFeraySniady2013} shows that
this is equivalent to the condition \ref{enum:types-of-edges}, 
as well as to the following condition:
\begin{multline*} |\F(M_{i})| - \text{(number of connected components of $M_i$)} = \\
   |\F(M_{i+1})| - \text{(number of connected components of $M_{i+1}$);}
\end{multline*}
in other words, 
the numbers 
\begin{equation}
\label{eq:difference}
 |\F(M_{i})| - \text{(number of connected components of $M_i$)}    
\end{equation}
over $\in\{0,1,\dots,n\}$ are all equal.
Since for $i=n$ the quantity \eqref{eq:difference} is equal to zero ($M_n=\emptyset$ is the empty map), 
this is equivalent to 
\begin{multline*}
 |\F(M_{i})| = \text{(number of connected components of $M_i$)}  \\  \text{for each } 1\leq i\leq n 
\end{multline*}
which is clearly the condition \ref{enum:maptopdegree}.
Remember that together with the removal of some leaf we always remove its endpoint, since we do not allow maps having isolated vertices, i.e.~connected components consisting of one vertex.

\smallskip
 
\emph{By now we have proved that the conditions \ref{cond:weight-max}, \ref{enum:types-of-edges}, 
and \ref{enum:maptopdegree} are all equivalent.}
On the other hand, condition \ref{enum:maptopdegree} implies that the map $M=M_0$ is top-degree
which is the condition \ref{cond:oneface}.
This concludes the proof of the equivalence.

\bigskip

Condition \ref{enum:types-of-edges} implies that each of the weights 
$\weight_{M_i,E_{i+1}}$ is a monic monomial in $\gamma$.
Condition \ref{enum:topdegree} implies that the product of these monic monomials
has the degree equal to $|\F(M)|+ |\E(M)| - |\V(M)|$.
This concludes the proof of \eqref{eq:monic}.
\end{proof}

\subsection{Measure of non-orientability of a map}
\label{subsec:measure-of-nonorientability}

Let $M$ be a map with $n$ edges. We define
\begin{equation}
\label{eq:take-average}
 \weight_M =\weight_M(\gamma) := \frac{1}{n!} \sum_{\prec} \weight_{M,\prec}.  
\end{equation}
This quantity can be interpreted as the mean value of the weight associated to
the map $M$ equipped with a randomly
selected history (with all histories having equal probability).
We call $\weight_M$ \emph{the measure of non-orientability} of the map $M$.

\begin{example}
\label{example:B}
We revisit \cref{example:A}.
For the histories $\{3,6\} \prec \{2,4\} \prec \{1,5\}$ 
and 
$\{3,6\} \prec \{1,5\} \prec \{2,4\}$ 
when the edge $\{3,6\}$ is removed first
the corresponding weight is equal to 
$\weight_{M,\prec}=1 \cdot \frac{1}{2} \cdot 1$.
For the remaining $4$ histories the corresponding weight is equal to 
$\weight_{M,\prec}=\gamma \cdot \gamma \cdot 1$. 
Finally, 
$$ \weight_M = \frac{2 \times 1 \cdot \frac{1}{2} \cdot 1+ 4 \times \gamma \cdot \gamma \cdot 1 }{6}.$$
\end{example}

\begin{proposition}
\label{prop:why-montop}
Let $M$ be a non-oriented map with $n$ edges.
The corresponding polynomial $\weight_M$ is of degree at most 
$n+|\F(M)|-|\V(M)|$.
The corresponding leading coefficient 
\[ \left[ \gamma^{n+|\F(M)|-|\V(M)|} \right] \weight_M = \weight_M^{\ttop}\]
is given by the quantity defined in \cref{def:mon-top-degree}.
\end{proposition}
\begin{proof}
This is an immediate application of \cref{lem:connected-component-one-face}.
\end{proof}

\subsection{Orientability generating series}

We recall that the orientability generating series was defined in \eqref{eq:what-is-ogs}
as a weighted sum of the normalized numbers of embeddings,
with the weight given by the polynomial $\weight_M$ described above:
\begin{equation}
 \widehat{\Ch}_\pi(\lambda) := (-1)^{\ell(\pi)} \sum_M \weight_M  \Embed_M(\lambda),
\end{equation}
where the sum is a conservative summation over non-oriented maps
with the face-type $\pi$.

\begin{lemma}
\label{lem:degree-bound}
For each partition $\pi$ there exists a Stanley polynomial for
$\widehat{\Ch}_\pi$; this Stanley polynomial is of degree $|\pi|+\ell(\pi)$.  
\end{lemma}
\begin{proof}
It is relatively easy to show that for any bicolored graph $G$
the corresponding normalized number of embeddings 
$\lambda\mapsto \Embed_G(\lambda)$ has a Stanley polynomial
which  is homogeneous of degree $|\V(G)|$.

\cref{prop:why-montop} provides an upper bound on the degree of the polynomial
$\weight_M$. 

By combining the above two observations it follows that
the Stanley polynomial for the product 
$\weight_M \Embed_M$ exists and its degree is bounded from above
by 
$ |\E(M)|+|\F(M)|= |\pi|+\ell(\pi),$
as required.
\end{proof}

\subsection{The top-degree of the orientability generating series}
\label{sec:newCh-as-a-polynomial}

\begin{corollary}
\label{coro:mon-top-degree}
For any integer $n\geq 1$ the top-degree homogeneous part of the orientability generating series is given by
\[   
\widehat{\Ch}^{\ttop}_n = \sum_M \weight_M^{\ttop} \gamma^{n+1-|\V(M)|}\ \Embed_M,
 \]
where the sum on the right-hand side runs over non-oriented maps with face-type $(n)$.
\end{corollary}
\begin{proof}
By \cref{lem:degree-bound} Stanley polynomial for $\widehat{\Ch}_n$ is of degree (at most) $n+1$;
out goal now is to extract its homogeneous part of degree $n+1$. 
This can be done by revisiting the proof of \cref{lem:degree-bound}
and using \cref{prop:why-montop}. 
\end{proof}

\subsection{Proof of \cref{theo:top-degre-equal}}
\label{sec:proof-second-main-result}

\begin{proof}[Proof of \cref{theo:top-degre-equal}]
The left-hand side of \eqref{eq:chtops-are-born-equal}
is given by  \cref{theo:formula-for-chtop}.
The right-hand side of \eqref{eq:chtops-are-born-equal}
is given by \cref{coro:mon-top-degree}.
Now it is enough to apply \cref{theo:main-bijection}
to show that they are equal.
\end{proof}

\section*{Acknowledgments}

We thank Maciej Dołęga and Valentin F\'eray for several years of collaboration on topics related to the current paper.

Research supported by \emph{Narodowe Centrum Nauki}, grant number \linebreak 2014/15/B/ST1/00064.

\bibliographystyle{alpha}
\bibliography{biblio}

\end{document}